\documentclass[12pt,fleqn,emtex]{article}

\usepackage{graphicx}
\usepackage{framed}
\usepackage{times}
\usepackage{mathptmx}

\usepackage[german, american]{babel}  

\usepackage{amsfonts}       %
\usepackage{amsmath} 
\usepackage{amssymb}        %
\usepackage{amstext}        %
\usepackage{amsthm}         
\usepackage[matrix,curve]{xypic}

\usepackage{epsfig}
\psfigdriver{dvips}
\usepackage{latexsym}
\usepackage[matrix,curve]{xypic}
\usepackage{rotating}
\rotdriver{dvips}

\newcommand{\Z}{\mathbb{Z}}
\newcommand{\Q}{\mathbb{Q}}
\newcommand{\R}{\mathbb{R}}
\newcommand{\C}{\mathbb{C}}

\newtheorem{Lemma}{Lemma}
\newtheorem{Theorem}{Theorem}
\newtheorem{Corollary}{Corollary}
\newtheorem{Proposition}{Proposition}

\DeclareMathSizes{12}{11}{8}{5} 
\usepackage{xypic,dsfont}

\begin{document}

\title{Degenerate Perverse Sheaves on Abelian Varieties}
\author{Rainer Weissauer}
\date{}

\maketitle
\thispagestyle{empty}

\centerline{\bf }

%

\bigskip
In [KrW] generic vanishing theorems for perverse sheaves on complex abelian varieties $X$ were proved. The proofs in loc. cit. depend on the characterization of perverse sheaves on $X$ with vanishing Euler-Poincare characteristics in terms of translation invariant perverse sheaves. We discuss this characterization in the present paper with its main result theorem 1 (for an overview see section 1). Although we deal with the situation in characteristic zero, 
at certain points we also use methods from characteristic $p>0$ also considered in [W3].

\bigskip
By sheaf theoretic methods, we show how the characterization of perverse sheaves on complex abelian varieties with vanishing  Euler-Poincare characteristic can be reduced to the study  of perverse sheaves on simple abelian varieties. This suffices for the proof of the vanishing theorems in [KrW], since the special case of simple complex abelian varieties has already been considered there. 

\bigskip
The proof for simple abelian varieties in loc. cit. uses the theory of $D$-modules, but this proof unfortunately does not carry over the case of non-simple abelian varieties and therefore has to be complemented by the arguments given here.  
We mention that independently another proof for the mentioned characterization was given by Chr. Schnell [S], also using the theory of $D$-modules.  The methods developed in the present paper use induction on the number of simple factor. They are not restricted to characteristic zero, but also work over finite fields -- so far with one 
exception,  the case of simple abelian varieties as the induction start. 
Hence, to generalize all results of [KrW] to abelian varieties over finite fields, it remains to characterize perverse sheaves with vanishing Euler-Poincare characteristic on simple abelian varieties over finite fields. 


\bigskip\noindent

\section{Formulation of the main results}\label{rel}

\bigskip\noindent

\medskip   
{\it Notations}. Let $X$ be a complex abelian variety. 
Let $E(X)$ denote the class of perverse sheaves whose irreducible constituents $K$ 
satisfy $\chi(K)=0$. Let $N(X)$ denote the class of {\it negligible} perverse sheaves on $X$, i.e. those for which each irreducible constituent is {\it translation invariant} for certain abelian subvarieties  of $X$ of dimension $>0$. A typical example for a translation invariant irreducible perverse sheaf on $X$ is $\delta_X^\psi = L_\psi[\dim(X)]$ for the local system $L_\psi$ on $X$
defined by a character $\psi: \pi_1(X,0)\to \C^*$ of the fundamental group of $X$.
Indeed, $T_x^*(\delta_X^\psi) \cong \delta_X^\psi$ holds for all $x\in X$, hence $\delta_X^\psi \in N(X)$. As already mentioned, it is not difficult to show $N(X) \subseteq E(X)$. If $F(X)$ denotes the class of irreducible perverse sheaves in the complement $E(X)\setminus N(X)$, then our aim is equivalent to show $F(X)=\emptyset$. Let $M(X)$ denote the class of  perverse sheaves whose irreducible components $M$ have Euler-Poincare characteristic $\chi(M)\neq 0$.   

\medskip
For a perverse sheaf $K$ on $X$
the character twist $K_\psi = K\otimes_{\C_X} L_\psi$ is a perverse sheaf. The classes
$N(X)$ and $E(X)$ are stable under character twists in this sense. Depending
on the situation, we also write $K^\psi$ instead of $K_\psi$, e.g.  $\delta_X^\psi=(\delta_X)_\psi$. 

\medskip 
For isogenies $f:X\to Y$ between complex abelian varieties the functors $f_*$ and $f^*$ preserve
$F(X)$, $E(X)$ and the categories of perverse sheaves; this is an easy consequence of the properties of the class $N_{Euler}$ (defined in [KrW]) and the adjunction formula. 
A complex $K$ with bounded constructible cohomology sheaves is called {\it negligible}, if its perverse cohomology sheaves ${}^p H^\nu(K)$ are all contained in $N(X)$.

\medskip
Our main result stated in theorem \ref{MAINTHEOREM}  is the assertion $F(X)=\emptyset $, respectively the equivalent

\begin{Theorem} \label{mn} For an irreducible perverse sheaf $K$ on 
a complex abelian variety $X$ with vanishing Euler characteristic 
$\chi(K) = \sum_i (-1)^i \dim(H^i(X,K))$ there
exists a nontrivial abelian subvariety $A\subseteq X$ such that 
$T_x^*(K) \cong K$ holds for all $x\in A$. 
\end{Theorem}

For simple complex abelian varieties this is shown in [KrW], using that
any perverse sheaf $K$ has an associated $D$-module whose characteristic variety as a subvariety of the cotangent bundle is a union of Lagrangians $\Lambda=\Lambda_Z$
for irreducible subvarieties $Z\subseteq X$. The assumption $\chi(K)=0$ implies
that all $Z$ are degenerate (see [W] and [KrW]), hence for simple $X$ this immediately  implies $Z=X$.
By the Lagrangian property then $\Lambda$ 
is the zero-section of the cotangent bundle $T^*(X)$. Hence by a well known theorem
on $D$-modules $K$ is attached to a local system defined and smooth on $X$, and then $K$ is a translation invariant perverse sheaf on $X$. This
proves the statement for simple abelian varieties, and this essentially is the proof of [KrW]. 
To give a similar argument for the non-simple case fails, or at least the author's attempt to give a simple proof along these lines for general abelian varieties. I would like to thank Christian Schnell for pointing this out.

\bigskip
In this paper we show how the general case of the theorem \ref{mn} can be reduced to the case of simple abelian varieties by sheaf theoretic arguments.  Since this may be of independent interest,
we mention that our arguments for the proof of theorem \ref{mn} carry over to abelian varieties and perverse sheaves over finite fields, once we assume that theorem 1 holds for simple abelian varieties and perverse sheaves over finite fields. In this case the Euler-Poincare characteristic 
has to be defined by extending $X$ and $K$ to the algebraic closure $k$ of the finite field. 
In fact, one step of our argument for characteristic zero (section 9 and appendix) even uses methods of characteristic $p$ by referring to [W3, appendix].

\medskip
{\it Some reformulations of theorem \ref{mn}}.
Simple perverse sheaves $K$ on $X$ are of the form $K = i_*(j_{!*}E_U)$
for some local system $E_U$ on an open dense subvariety $j: U \hookrightarrow Z$ 
of the support $Z=supp(K)$. The support is an irreducible closed 
subvariety $i: Z \hookrightarrow X$
of $X$. If the irreducible perverse sheaf $K$ is translation invariant with respect to an abelian
subvariety $A$ in $X$, its support $Z$ satisfies $Z+A=Z$.
The local system $E_U$ defines an irreducible finite dimensional representation $\phi$ of the fundamental
group $\pi_1(U)$ of $U$.  By the $A$-invariance of the singular support of $K$
there exists an $A$ invariant closed subset $Z'$ of $Z$, which contains the ramification locus of the perverse sheaf $K$. The restriction of the perverse $K$ to $U=Z\setminus Z'$ is smooth in the sense that $K\vert_U = {E} [d]$ holds for a smooth etale sheaf $E$ associated to the representation of the topological fundamental group $\pi_1(U)$ of $U$. Notice, $U$ can be chosen such that $U
+A=U$ holds.
 
 \medskip
 Let $\tilde U$ and $\tilde Z$ denote the images of $U$ and $Z$ under the projection  $q: X\to \tilde X=X/A$.
Since $q: X \to \tilde X=X/A$ is smooth, also the induced morphism
$q : Z \to \tilde Z$ is smooth. 
The smooth morphism $q:Z \to \tilde Z$ defines a fibration $q: U\to \tilde U$. Since $A$ is connected, we obtain from the long exact homotopy sequence
$$ \xymatrix{  \pi_1(A) \ \ar[r]^\sigma & \ \pi_1(U)\  \ar[r] & \ \pi_1(\tilde U)\ \ar[r] & 0 }\ .$$ 
The map $\sigma: \pi_1(A)\to \pi_1(U)$ is injective. To show this, consider the natural group homomorphism $\rho: \pi_1(U) \to \pi_1(X)$
induced from the inclusion $U \hookrightarrow X$. Obviously the composition $\rho\circ \sigma$ 
is the first map of the exact homology sequence
$$ \xymatrix{ 0 \ \ar[r] &\ \pi_1(A) \ \ar[r] & \ \pi_1(X)\  \ar[r] & \ \pi_1(\tilde X)\ \ar[r] & 0 }\ .$$ 
Hence $\rho\circ \sigma$ (and therefore $\sigma$) is injective.
$\pi_1(A)$ is a normal subgroup of $\pi_1(U)$.
We claim that $\pi_1(A)$ is in the center of $\pi_1(U)$. Indeed, for $\alpha\in \pi_1(A)$ and $\gamma\in \pi_1(U)$, there exists an $\alpha'\in \pi_1(A)$ such that
 $\gamma \alpha \gamma^{-1} = \alpha'$. If we apply the homomorphism $\rho$,
this gives $\rho(\gamma) \rho(\alpha) \rho(\gamma)^{-1} = \rho(\alpha')$. Hence $\rho(\alpha) = \rho(\alpha')$, since
$\pi_1(X)=H_1(X)$ is abelian. Therefore $\alpha=\alpha'$, because $\rho\circ \sigma$ is injective.
Because $\pi_1(A)$ is a central subgroup of $\pi_1(U)$, for any irreducible representation $\phi$ of $\pi_1(U)$ there exists a character $\chi$ of $\pi_1(A)$ such that $\phi(\alpha\gamma)=\chi(\alpha)\phi(\gamma)$ holds for $\alpha\in \pi_1(A)$ and $\gamma\in \pi_1(U)$. Since $\pi_1(\tilde X)$ is a free $\Z$-module, any character $\chi$ of $\pi_1(A)$ can be extended to a character $\chi_X$ of $\pi_1(X)$.
Thus $\chi_X^{-1} \otimes\phi$ is an irreducible representation, which is trivial on $\pi_1(U)$; in other words it is an irreducible representation of $\tilde U$.

\medskip
The last arguments imply that there exists a perverse sheaf $\tilde K$ on $\tilde U$
such that $L_{\chi_X}^{-1}\otimes K = q^*(\tilde K)[\dim(A)]$ holds on $U$. Then
$\tilde K$ necessarily is an irreducible perverse sheaf. 
Let $\tilde K$ also denote the intermediate extension of $\tilde K$ to $\tilde Z$, which is an irreducible perverse sheaf on $\tilde Z$. Since $q:Z\to \tilde Z$ is a smooth morphism with connected fibers, the 
pullback $q^*[dim(A)]$ is a fully faithful functor from the category of perverse sheaves on $\tilde Z$ to the category of perverse sheaves on $Z$. Also $L=q^*(\tilde K)[\dim(A)]$, as a perverse sheaf on $Z$, is still irreducible on $Z$. Now $K$ and $L$ are both irreducible perverse sheaves on $Z$, whose restrictions on $U$ coincide. Thus $K=L$.  
Choose a finite etale covering
such that $\tilde X$ splits. Then this implies the assertion (b) $\Longrightarrow$ (c) of the next theorem. The implications (b) $\Longrightarrow$ (c) $\Longrightarrow$ (a) of the next theorem are elementary. Hence, we have shown that quite generally theorem \ref{mn} implies 

\medskip
\begin{Theorem}  For an irreducible perverse sheaf $K$ on a complex abelian variety $X$ the following
properties are equivalent

\medskip
\begin{itemize}
\item[a)] The Euler characteristic $\chi(K)$ vanishes.

\smallskip
\item[b)] There exists a positive dimensional abelian subvariety $A$ of $X$, a translation invariant smooth sheaf $L_{\chi_X}$ of rank one on $X$ and
and a perverse sheaf $\tilde K$ on $\tilde X=X/A$ such that $K \cong L_{\chi_X} \otimes q^*(\tilde K)[\dim(A)]$ holds for the quotient map $q:X\to \tilde X$.

\smallskip
\item[c)] There exists a finite etale covering of $X$ which splits into a product of two abelian varieties $A$ and $\tilde X$, where $\dim(A)>0$, such that  the pullback of $K$ to this covering
is isomorphic to the external tensor product of a translation invariant perverse sheaf on $A$ and a perverse
sheaf on $\tilde X$.
\end{itemize}
\end{Theorem}

\medskip 
{\it Outline of the proof of theorem \ref{mn}}.
The proof is by induction on the number of simple factors of $X$ (up to isogeny). 
The case of simple abelian varieties we now take for granted. 
For the proof we use that the perverse sheaves $K$ may be replaced by monoidal perverse sheaves on $X$ in the sense of [W2]. 
For $X$ isogenous to a product $A_1\times A_2$ of two simple abelian varieties  we use methods from characteristic $p$ in the proof. We deal with this case in section \ref{specialsheaves} after some preparations in section \ref{previous} and \ref{next} relying on arguments that involve the tensor categories introduced in [KrW]. 

\medskip
If $X$ has three or more  simple factors, we will simply use
induction on $dim(X)$. The main step in this case, obtained in section \ref{main}, uses an analysis of the stalks of monoidal perverse sheaves. The arguments are sheaf theoretic, exploiting
spectral sequences that naturally arise if one restricts perverse sheaves on $X$ to abelian subvarieties (e.g. fibers of homomorphisms) used in the course of the induction argument. These spectral sequences are described in section \ref{steps} and are applied in section \ref{FS} and \ref{main}. A second crucial step consists in the  use of monoidal components ${\cal P}_K$ of perverse sheaves $K$ on $X$ and their degree $\nu_K$ (see [W2]). Furthermore,
in section \ref{FS} we show that it suffices to consider perverse sheaves in $F_{max}(X) \subseteq F(X)$. 

\bigskip\noindent

\section{Relative generic vanishing}\label{rel}

\bigskip\noindent

\medskip
An important technical tool for the study 
of homomorphisms $f:X\to B$ used through the induction process will be the next lemma. This lemma can be easily derived from the statement $F(A)=\emptyset$ for simple abelian varieties $A$, whose proof 
for complex abelian varieties was sketched above.
The statement $F(A)=\emptyset$ can be converted into a weak relative generic vanishing theorem for morphisms with simple kernel $A$ (see [KrW]). 
Factoring arbitrary homomorphisms  $f:X\to B$ into homomorphisms whose kernels are simple abelian varieties, an iterative application of the assertion  
$F(A)=\emptyset $ for each of the simple abelian varieties $A$ defining the sucessive kernels, then easily implies
the next

\goodbreak

\begin{Lemma} \label{relGVT} Let $A$ be an abelian subvariety of $X$ with $F(A)=\emptyset$ and quotient map
$$ f: X \to B=X/A \ $$
and let $K\in Perv(X,\C)$ be a perverse sheaf on $X$. Then for a generic character $\chi$  the direct image $Rf_*(K_\chi)$ is a perverse sheaf on $B$.
\end{Lemma}

\medskip
The similar statement for abelian varieties over finite fields follows from the vanishing theorem proved in [W4]. 

\medskip
Under the assumptions of the last lemma, then $K\in E(X)$ immediately implies 
$Rf_*(K_\chi)\in E(X/A)$. Inductively this gives 

\begin{Corollary}\label{azyk}
$H^\bullet(X,K_\chi)=0$ holds for $K\in E(X)$ and generic $\chi$.
\end{Corollary}

\medskip
An irreducible perverse sheaf $K$ on $X$ is called  {\it maximal}, if for any 
quotient homomorphism $f:X\to B$ to
a simple abelian quotient variety  $B$ and a generic character twists $K_\chi$ of $K$ the direct image $Rf_*(K_\chi)$ does not vanish. Let $F_{max}(X)$ denote the maximal perverse sheaves 
in $F(X)$.  If , then for $K\in F_{max}(X)$ one has
the following assertion. 

\begin{Lemma} Suppose $K\in F_{max}(X)$ and suppose given a surjective homomorphism $f:X\to B$ for simple $B$. If $F(B)=\emptyset $, then $Rf_*(K_{\chi_0})\neq 0$ holds for any character $\chi_0$. 
\end{Lemma}

\medskip
{\it Proof}. By corollary \ref{azyk}, for generic $\chi_0$ the direct image complex 
$L=Rf_*(K_{\chi_0})\neq 0$ is a translation invariant perverse sheaf on $B$. Therefore its cohomology sheaves are concentrated in degree $-\dim(B)$. So, the same holds for the cohomology of $M:= K_{\chi_0}\vert_{f^{-1}(b)}$ over the fiber
$f^{-1}(b)$ for any $b\in B$. Hence the Euler-Poincare characteristic of $M$ does not vanish.
Since the Euler-Poincare characteristic is invariant under characater twists, the same same holds for any character $\chi_0$. This implies the assertion of the lemma. \qed

\medskip
{\bf Remark}. \label{pseudo} From now on we use that $F(A)=\emptyset$ holds for all simple complex abelian varieties. For finitely many perverse sheaves on $X$  and a given homomorphism $f:X\to B$ to an arbitrary abelian variety $B$, then 
 lemma \ref{relGVT} allows to find characters $\chi$ for which  $\Gamma^\chi(K)
= {}^p H^0(Rf_*(K_\chi))$ defines an exact functor $\Gamma^\chi$ on the tensor subcategory ${\cal T}$ of $D_c^b(X,\C)$ generated by these objects with respect to the  convolution product (see [KrW]) , in the sense that $\Gamma^\chi$  maps distinguished triangles
to short exact sequences. We can use this to analyse stalks:
Suppose the stalk $Rf_*(K_\chi)_b$ vanishes for generic $\chi$. Let $F=F(b)$
denote the fiber $f^{-1}(b)$. Then the restriction $M := (K_\chi)\vert_{F(b)}$ is in ${}^p D^{[-\dim(B),0]}(F)$ and
for generic $\chi$ all its perverse cohomology sheaves $M^i = {}^p H^{-i}(M)$ are acyclic, i.e. $H^\bullet(F,M^i)=0$. Although the sheaf complex $M$ and also the perverse sheaves $M^i$ are not necessarily semisimple, this assertion easily follows from the exactness of the functors 
$\Gamma^\chi$, here as a consequence of
$$ H^i(F,M) \ = \ H^0(F,M^i)\ =\ H^\bullet(F,M^i) \ $$
for generic $\chi$ (identifying $F$ and $A=kern(f)$). 
The same statement carries over to the finitely many irreducible perverse Jordan-H\"older constituents $P$ of the perverse sheaves $M^i$. 

\medskip
{\it Direct images}. For the definition of ${M}$ above we fixed a suitable generic character $\chi$. Then  ${M}_{\chi_0} = K_{\chi\chi_0}\vert_{F(b)}$ (for arbitrary $\chi_0$) gives the \lq{\it cohomology}\rq\ spectral sequence 
$$ \bigoplus_{j+ i = k} H^{-j}(F(b), (M^i)_{\chi_0}) \Longrightarrow
R^{-k}f_*(K_{\chi\chi_0})_b \ .$$
The cohomology sheaves
${\cal H}^{i+k}(M^i) =  
 \bigoplus_i {\cal H}^{k}(M^i[i])$ are related to the 
 cohomology sheaves ${\cal H}^k({M})$, or equivalently the stalk cohomology sheaves of the complex $K_\chi$ at points $x\in F(b)$ via
 the following \lq{\it stalk}\rq\ spectral sequence with $E_2^{p,q}= {\cal H}^p(M^{-q}) $
$$   \bigoplus_{-p-q=-k} {\cal H}^{-p}(M^{q}) \ \ \Longrightarrow \ \
 {\cal H}^{-k}({M}) = {\cal H}^{-k}(K_\chi)\vert_{F(b)} \ \ $$
on $F(b)$ with differentials $ d_2:  {\cal H}^{i+ k-1}(M^{i-1})  \to  {\cal H}^{i+k+1}(M^i)$.

$$\xymatrix@-0,9cm{
{\cal H}^{-d(A)}(M^0)\ \  ... &  {\cal H}^{-d-1}(M^0) & {\cal H}^{-d}(M^0) &   &  {\cal H}^{-1}(M^0)\ar@{-}[ll]   &  {\cal H}^{0}(M^0)\ar@{-}[l] & 0\ \ ... \cr
{\cal H}^{-d(A)}(M^1)\ \  .... & {\cal H}^{-d-1}(M^1) & {\cal H}^{-d}(M^1) & ...   &  {\cal H}^{-1}(M^1)   &  {\cal H}^{0}(M^1) \ar@{-}[u]& 0\ \ ... \cr
   \ \ ...    &       .          &    &    &  . &                         & 0 \ \ ...     \cr
    \ \ ... &           .            &  &    &  . &                       & 0 \ \ ...       \cr
  \ \ ... &             .            &   &    &  . &                       & 0 \ \ ...        \cr
{\cal H}^{-d(A)}(M^d) \ \ ... & {\cal H}^{-d-1}(M^d) &  {\cal H}^{-d}(M^d) & ...   &  {\cal H}^{-1}(M^d) &  {\cal H}^{0}(M^d) \ar@{-}[llluuuuu]\ar@{-}[uuuu] & 0 \ \ ... \cr
 \ \ ... &             .            &  . &  .  &  . &               .        & 0 \ \ ...        \cr
  \ \ ... &             .            &  . &  .  &  . &               .        & 0 \ \ ...        \cr
{\cal H}^{-d(A)}(M^{d(B)}) \ \ ... & {\cal H}^{-d-1}(M^{d(B)}) &  {\cal H}^{-d}(M^{d(B)}) & ...   &  {\cal H}^{-1}(M^{d(B)}) &  {\cal H}^{0}(M^{d(B)})  & 0 \ \ ... 
} $$

\medskip\noindent
Here $d(A)$ and $d(B)$ are the dimensions of $A$ resp. $B$. The 
edge morphisms ${\cal H}^0(M^i)\to {\cal H}^{-i}({M})$ and
${\cal H}^{-i}({M}) \to {\cal H}^{-i} (M^0)$ will play a role later, as well as the upper triangle of the diagram
(as indicated) for a suitable $d$.

\begin{Lemma} \label{notzero}
For simple $A\subseteq X$ and $f:X\to X/A$ suppose $F(A)=\emptyset$. Then for an irreducible  perverse sheaf $K$ on $X$ we have $ Rf_*(K)=0\ \ \Longrightarrow \ \ K\notin F(X)$.
\end{Lemma}

{\it Proof}. $Rf_*(K)=0$ implies $Rf_*(K_\chi)=0$ for generic $\chi$ (proof of corollary \ref{azyk}). So all $M^i$  are acyclic perverse sheaves on $A$, hence contained in $E(A)$.
Now use $F(A)=\emptyset$. Since $A$ is a simple abelian variety $A$, it implies  ${\cal H}^{-j}(M^i)=0$ for $j\neq d(A)$ by translation invariance. Thus
 by the stalk spectral sequence, ${\cal H}^{-i-d(A)}(K_\chi)\vert_{F(b)} \cong {\cal H}^{-d(A)}(M^i)$ is
a translation invariant sheaf on $A$ and therefore is never a skyscraper sheaf.
We now apply this for the monoidal component ${\cal P}_K$ of $K$. According to [W2, lemma 2.1],
$Rf_*(K_\chi)=0$ for generic $\chi$ implies $Rf_*({\cal P}_{K\chi})=0$ for generic $\chi$ and  $K\in F(X)$ implies ${\cal P}_K \in F(X)$.
For ${\cal P}_K \in F(X)$  the stalk ${\cal H}^{-i}({\cal P}_K)$ is a skyscraper sheaf at least for the degree $i=\nu_K$ by [W2, lemma 1, part 7]. This gives a contradiction if $K\in F(X)$.
\qed


%

\medskip

\section {Restriction in steps}\label{steps}
 
In this section we consider exact sequences of abelian varieties
$$ \xymatrix@R-1cm{ 0 \ar[r] & B_1 \ar[r] & B \ar[r]^h & B_2 \ar[r] & 0} \ $$
$$ \xymatrix@R-1cm{ 0 \ar[r] & C \ar[r] & A \ar[r]^{p} & B_1 \ar[r] & 0} \ $$  
and a diagram of quotient homomorphisms, where $A=g^{-1}(B_1)$
and $p=g\vert_A$ and where $C$ is the kernel of the projection $g:X \to B$ 
$$  \xymatrix@C+0,3cm{  X \ar[rr]^f \ar[dr]_g & & \ \ B_2=X/A \cr
& B=X/C \ar[ur]_h & } $$
 

\medskip
{\bf Assumption}. \label{ass} {\it For 
perverse $K$ and
a quotient homomorphism $f: X \to B_2=X/A $ 
 suppose for generic $\chi$ and some $d$ } 
$$ \fbox{$ R^{-i}f_*(K_\chi) \ =\  0 \quad , \quad  \forall \ i < d $} \ .$$
These vanishing conditions imply acyclicity for the constituents $P$ of the perverse sheaves $M^0,M^1,...,M^{d-1}$.


\medskip
For $b_2 \in B_2 $ and $b_1 \in B_1$ the fibers
$  C \cong F(b_1,b_2) = g^{-1}(b_1,b_2)  \hookrightarrow    f^{-1}(b_2) = F(b_2) \cong A  $,
can be identified with $A$ respectively with
$C$ up to a translation. This being said, we restrict a generic twist $K_\chi$ of $K$ (for some {generic} character  $\chi: \pi_1(X,0) \to \C^* $) to
the fiber $F_{b_2}$ 
$$  {M} = {M}(b_2) = (K_\chi)\vert_{F(b_2)}   \quad \in \ \ {}^p D^{[-\dim(B_2),0]}(F(b_2)) \ ; $$
then we further restrict ${M}$  to $  F(b_1,b_2)  \hookrightarrow  F(b_2) $ and obtain
$$  {N} = {N}(b_1,b_2) := {M}(b_2)\vert_{F(b_1,b_2)} = K_\chi\vert_{F(b_1,b_2)} \ \in {}^pD^{[-\dim(B),0]}(F(b_1,b_2))\ .$$ 
For $N^k= N^k(b_1,b_2) = {}^pH^{-k}({N})$ in $Perv(F(b_1,b_2),\C)$ and $M^i = M^i(b_2)= {}^p H^{-i}({M})$
in  $Perv(F(b_2),\C)$
for  $j=0,..,\dim(B_1)$ and $i=0,..,\dim(B_2)$ we easily get the following \lq{\it double restriction}\rq\ spectral sequence 
$$   \bigoplus_{i+j = k} {}^p H^{-j}\bigl(M^i(b_2)\vert_{F(b_1,b_2)} \bigr) 
\Longrightarrow  N^k(b_1,b_2)  \ $$


\medskip
{\it Picture}. The front rectangle visualizes the fiber $F(b_2)$, which is isomorphic to $A=Kern(f)$,
and the fiber $F(b_1,b_2) \subset F(b_2)$  isomorphic to the abelian variety $C=Kern(g)$. 
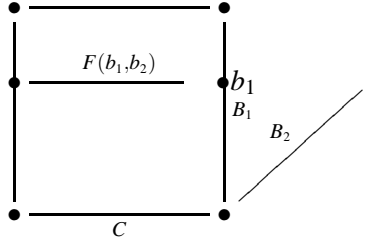
$$ \xymatrix@-0,3cm{ \bullet & & & \bullet \ar@{-}[lll] & & \cr
\bullet  & & & \ \ \ \bullet b_1 \ar@{-}[lll]_{F(b_1,b_2)} & & \cr
 & & & & & \cr
 \bullet \ar@{-}[uuu] & & & \bullet \ar@{-}[lll]^C \ar@{-}[uuu]_{B_1}\ar@{-}[uurr]^{B_2} & &
}$$

\medskip\noindent
Now fix $b_2\in B_2$. Then for almost all closed points $b_1\in B_1$ the perverse sheaf $${}^p H^{0}(M^d(b_2)\vert_{F(b_1,b_2)})$$ is zero, since it defines a perverse
quotient sheaf of $M^d(b_2)$ on $F(b_2)$ with support in $F(b_1,b_2)$.  Indeed these supports are disjoint and there are only finitely many constituents. Furthermore the perverse constituents of the sheaves $M^0(b_2),..,M^{d-1}(b_2)$
on $F(b_2)\cong A$ are in $E(F(b_2))$, by the vanishing assumption on the direct images: $Rf^{-i}(K_\chi)=0$ for $i<d$. 

\medskip
For generic $\chi$ we have the 
\lq{\it relative cohomological}\rq\ spectral sequence
$$ \bigoplus_{j+ i = k} H^{0}(F(b_1,b_2), {}^pH^{-j}(M^i\vert_{F(b_1,b_2)}) )\Longrightarrow 
R^{-k}g_*(K_{\chi})_{(b_1,b_2)} \ .$$
Notice $i=0,1,...,\dim(B_2)$ and $j=0,...,\dim(B_1)$,
where the case $j=0$ plays a special role as explained above.
The spectral sequence is obtained from the double restriction spectral sequence
combined with the degenerate cohomology spectral sequence for $g$
using $H^{0}(F(b_1,b_2), N^k) = 
R^{-k}g_*(K_{\chi})_{(b_1,b_2)} $ for generic $\chi$.

\medskip  
Now assume $$ \fbox{$ d=\dim(B_1)=\mu(A)= \mu(X) $} \ .$$

\begin{Proposition}\label{gut} For $d=\mu(A)=\dim(B_1)=\mu(X)$  suppose  given an irreducible perverse sheaf $K$ with the vanishing condition $R^{-i}f_*(K_\chi)_{b_2}=0$ for $i<d$ and generic $\chi$.
Then for fixed $b_2\in B_2$ and generic $\chi:\pi_1(X,0)\to \C^*$ we have an exact sequence
$$  H^{0}\bigl(F(b_1,b_2), {}^pH^{0}(M^d\vert_{F(b_1,b_2)}) \bigr) \to R^{-d}g_*(K_{\chi})_{(b_1,b_2)} \to H^0\bigl(F(b_1,b_2), {}^p H^{-d}(M^0\vert_{F(b_1,b_2)})\bigr) \ .$$
In particular $R^{-d}g_*(K_{\chi})_{(b_1,b_2)} =  0 $
holds for almost all $b_1\in B_1$ (for fixed $b_2\in B_2$),
if $M^0$ vanishes. Notice $M^0=0$ iff $K$ does not have support in the fiber
$F(b_2)$. 
\end{Proposition}

{\it Proof}. $\bigoplus_{j+ i = d} H^{0}(F(b_1,b_2), {}^pH^{-j}(M^i\vert_{F(b_1,b_2)}) )\Longrightarrow
R^{-d}g_*(K_{\chi})_{(b_1,b_2)}$ for generic $\chi$ and $k=d$ 
degenerates by Lemma 
\ref{L2}  which shows that for $0<i< d$ we can ignore all terms $   j =1,...,d-1 = \dim(B_1)-1 $ in this spectral sequence.
Since $j+i=k$,
for $k=d$ only the terms $(j,i)=(0,d)$ and the 
term $(j,i)=(d,0)$ remain. This proves our assertion. \qed

\medskip
Before we give the proof of the lemma, recall
that the abelian variety $A$ can be identified with the \lq{front rectangle}\rq\ $F(b_2)$, the fiber of $b_2$ for fixed $b_2\in B_2$, which contains $F(b_1,b_2)$
$$ 
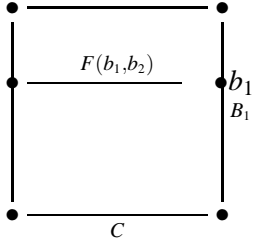
\xymatrix@-0,3cm{ \bullet & & & \bullet \ar@{-}[lll] & & \cr
\bullet  & & & \ \ \ \bullet b_1 \ar@{-}[lll]_{F(b_1,b_2)} & & \cr
 & & & & & \cr
 \bullet \ar@{-}[uuu] & & & \bullet \ar@{-}[lll]^C \ar@{-}[uuu]_{B_1} & &
}$$
The irreducible constituents $P$ 
of the perverse sheaves $M^i,i<d$  are acyclic perverse sheaves
living on the \lq{front rectangle}\rq\ $A$ and their irreducible perverse constituents $P$ are acyclic. 

\begin{Lemma} \label{L2} Suppose $\mu(A)=\dim(B_1)=\mu(X)=d$. Then $B_1$ is simple and for the constituents $P$ of $M^i$ for $i=0,1,..,d-1$ and  
for generic  $\chi$ the following holds: 
$$   H^{\bullet}\bigl(F(b_1,b_2),{}^p H^{-j}(P\vert_{F(b_1,b_2)})_\chi\bigr) \ = \ 0  \quad , \quad  \forall\ j=0,1,..,d-1 \ .$$
\end{Lemma}

{\it Proof}. The irreducible constituents $P=\tilde P_\chi$ (for $K\vert_{F(b_2)}$) of the $M^i$ for $i=0,1,..,d-1$ are acyclic on $A$ and for  $p:A \to B_1$ the direct image $Rp_*(P)$ is perverse ($\chi$ being generic) so that therefore $\chi(Rp_*(P))=\chi(P)=0$ holds. Since $B_1$ is simple of dimension $\dim(B_1)=d$, the semisimple perverse sheaf $Rp_*(P)$ in $E(B)$ is of the form  $$Rp_*(P)=\bigoplus_{\psi\in \Psi(\chi)} m_\psi\cdot \delta_{B_1}^\psi\ .$$ Then $R^{-i}p_*(P)_{b_1}=0$ for $i\neq d=\dim(B_1)$ and all $b_1\in B_1$, so that for generic $\chi$
the perverse sheaves  ${}^p H^{-j}(P\vert_{F(b_1,b_2)})$ are acyclic for $j=0,...,-d+1$ by the remark on page \pageref{pseudo} applied for $(P,p)$ instead of $(K,f)$. Thus $H^{\bullet}\bigl(F(b_1,b_2),{}^p H^{-j}(P\vert_{F(b_1,b_2)})_\chi\bigr)$ vanishes for $j=0,..,d-1$. \qed

\section{Stalk vanishing conditions}

Let $\mu(X)$ be the minimum 
of the dimensions of the simple abelian variety quotients $B$ of $X$. An 
abelian quotient variety  $B$ of $X$ is said to be {\it minimal} if $\dim(B)=\mu(X)$.  
For a sheaf complex $P$ on $X$ define 
$$    \mu(P) \ =  \ \max\{\nu \ \vert \ {\cal H}^{-i}(P) = 0 \ \mbox{ for all } \ i < \nu  \} \ .$$

%

\begin{Lemma}\label{A} $\mu(K) \geq \mu(X)$ holds for all perverse sheaves $K\in E(X)$. 
\end{Lemma}

\medskip
{\it Proof by induction on $dim(X)$}.
 Choose $f:X \to B$ with simple minimal $B$.
Then $Rf_*(K_\chi)$ is perverse for generic $\chi$, hence $Rf_*(K_\chi)\in E(B)$ is 
of the form $\bigoplus_\psi  m_\psi\cdot \delta_B^\psi$.
By the induction hypothesis we can assume $M^0=0$,
since otherwise the support of $K$ is contained in a proper abelian subvariety of $X$. 
${M}=K_\chi\vert_{f^{-1}(b)}$ has acyclic perverse cohomology
sheaves $M^i \in E(Kern(f))$ for $i=1,...,d-1$ and $d=\mu(X)$.
Then, by induction 
$\mu(M^i) \geq \mu(Kern(f))\geq \mu(X)=d$ implies
$   {\cal H}^{-\nu}(M^i) =0 $ for all $ i=0,...,\dim(B) -1=d-1$ and all $ \nu < d=\mu(X)$.
Hence $\mu(K) \geq d$ by the stalk
spectral sequence discussed in section \ref{rel}.
 \qed 

%
%
%

\medskip
To an irreducible perverse sheaf $K$ on $X$ we associated in [W2] a degree $\nu_K$ and an irreducible monoidal
perverse sheaf ${\cal P}_K$ with more amenable properties than $K$. It is defined from the monoidal structure obtained from the convolution product $*$ that comes from the group law on $X$. Then the complex ${\cal P}_K[-\nu_K]$ is the unique irreducible summand of
the convolution product $K^\vee*K$ on which the evaluation morphism is nontrivial (for the rigid dual $K^\vee$ of $K$).
In the next lemma we gather information on these monoidal components ${\cal P}_K$ for $K\in F(X)$. For details on monoidal components we refer to [W2].

\begin{Lemma} \label{Kill1} For $K\in F(X)$ the following holds
\begin{enumerate}
\item The monoidal
component ${\cal P}(K)$ is in $F(X)$ with $\nu_K:=\mu({\cal P}_K) \geq\mu(X)$. 
\item $K \in F_{max}$ implies ${\cal P}_K\in F_{max}$ 
and $\nu_K = \mu({\cal P}_K) = \mu(X)$.
\item For $K\in F_{max}(X)$  and any minimal quotient $B=X/A$ of $X$ 
the restriction ${\bf M}= {\cal P}_K\vert _F$  
of ${\cal P}_{K}$ to any fiber $F=f^{-1}(b)$ 
is a complex with Euler perverse cohomology $M^i$
for $i=0,..,\mu(X)-1$. For the fiber over the point $b=0$ 
furthermore $M^{\mu(X)}$
has a nontrivial perverse skyscraper quotient 
concentrated at the point zero.
\end{enumerate}
\end{Lemma}

\medskip
{\it Proof}. Concerning the first assertion ${\cal P}(K)\in F(X)$ see [W2, lemma 2.5].
Next, $\mu({\cal P}_K) \geq \mu(X)$ holds by lemma \ref{A}, by
$\nu_K := \mu({\cal P}_K)$ and [W2, lemma 1, part 7]. Hence $\nu_K\geq \mu(X)$.
By [W2, lemma 4], on the other hand $\nu_K\leq \mu(X)$ 
for $K\in F_{max}(X)$. Hence $\nu_K=\mu({\cal P}_K) = \mu(X)$ holds for $K\in F_{max}(X)$.
Since for ${\cal P}_K\in F(X)$ the cohomology ${\cal H}^{-\nu(K)}({\cal P}_K)$
is a skyscraper sheaf, our assertion on the skyscraper subsheaf comes 
from considering the edge term of the above spectral sequence. 
Since  $Rf_*(K_\chi)$ 
is perverse for generic $\chi$ and hence
$Rf_*(K_\chi)=0$ iff $Rf_*({\cal P}_{K\chi})=0$  by [W2, lemma
2.1], $K\in F_{max}(X)$ implies ${\cal P}_K\in F_{max}(X)$
\qed

\begin{Lemma} \label{maximal-}
For $K\in F_{max}(X)$ the support of $K$ is not contained in a translate of a proper abelian subvariety $A$ of $X$.
\end{Lemma}
 
{\it Proof}. 
For the projection $f:X\to X/A$ 
the support $Z$ of $Rf_*({\cal P}_K)$ becomes zero: $f(Z)=\{0\}$. 
This also holds for generic twists of $K$, 
so we could assume that $Rf_*({\cal P}_K)$ is perverse. Therefore
$Rf_*({\cal P}_K)=0$, because otherwise for a skyscraper sheaf 
$\chi(Rf_*{\cal P}_K)>0$ would hold; a contradiction. This implies
$Rf_*(K_\chi)=0$ for generic $\chi$, contradicting the maximality of $K$. \qed

\medskip

\section{Supports}\label{FS}
 
\medskip
%
%

Let $A$ be an abelian subvariety of $X$ 
and $K$ be an irreducible perverse sheaf on $X$. 
For quotient homomorphisms
$ f : X \to B=X/A$ consider the assertions
\begin{enumerate}
\item $K$ is $C$-invariant for some nontrivial abelian subvariety $C$ of $A$.
\item $Rf_*(K_\chi)=0$ for generic $\chi$.
\item $Rf_*(K_\chi)=0$ on a fixed dense open 
subset $W$ of $f(Z)$ for generic $\chi$, where $Z$ denotes the support of $K$.
\end{enumerate}
Obviously $1.\Longrightarrow  2.\Longrightarrow 3$. We remark that $Rf_*(K_{\chi_0})=0$ for a single $\chi_0$
implies 2. using corollary \ref{azyk}. Indeed the fibers ${M}$ are acyclic
and remain acyclic for generic character twist.

\begin{Proposition} \label{fantom}
Suppose $F(A)=\emptyset$ and $A+Z=Z$ for the support $Z$ of 
$K$. Then the three properties 1,2 and 3 above are equivalent.
Furthermore, for $K\in F(X)$ then  $supp(Rf_*(K_\chi)) = f(supp(K))$ holds for generic $\chi$. 
\end{Proposition}

\medskip
Step 1) Fix a smooth dense open subset $W \subset f(Z)=Z/A$ of dimension $d$ so that $K\vert_U = E[\dim(U)]$ for a local system $E$ on a dense Zariski open subset $U$ of $f^{-1}(W)$, and so that $U\cap f^{-1}(b)$ is dense in every fiber for $b\in W$; fix a closed point $b$ in $W$. Then $F=f^{-1}(b)$ can be identified with $A$.  
For $  {M}= (K_\chi)\vert_{F} \in {}^p D^{[-d,0]}(F)$ and $M^i = {}^p H^{-i}(M)$ for $i=0,..,d$ notice $  M^d \neq 0 $,
since ${}^p H^{-d}(M)$  contains the intermediate extension of $E_{F\cap U}[g] $ as constituent.

\medskip
Step 2) Let $C$ be a nontrivial abelian subvariety of $A$
and $g: X\to X/C$ the projection. By assumption, $K$ is simple and $g$ is smooth with connected fibers of dimension $\dim(C)$, so by [BBD, p.108ff] 
$K$ or equivalently $K_{\chi}$ is $C$-invariant for some nontrivial abelian subvariety $C$ of $A$ and some character $\chi$
iff
$$ {}^p H^{-\dim(C)}(Rg_*(K_\chi)) \neq 0 \ .$$

\medskip
Step 3) For fiber inclusions $i:F\hookrightarrow Z$ and $i_C: F/C \hookrightarrow Z/C$
$$ \xymatrix@+0,2cm{ Z \ar@/^10mm/[rr]^f\ar[r]^g &  Z/C \ar[r] & Z/A \cr
F \ar[r]^g\ar@{^{(}->}[u]^i & F/C \ar@{^{(}->}[u]^{i_C} \ar[r] & \{b\} \ar@{^{(}->}[u]^{i_b} } $$
proper base change gives $i_C^*Rg_*(K_\chi)
= Rg_*(i^*(K_\chi)) = Rg_*({M})$. In order to compute ${}^p H^{-\dim(C)-d}(Rg_* i^*(-))$, we use that under the functor $Rg_*$ perversity drops  at most by $-\dim(C)$ (see [BBD, 4.2.4])  $$Rg_*:\  {}^p D^{\geq n}(-)\ \longrightarrow \ {}^p D^{\geq n-\dim(C)}(-)$$ 
and that therefore
$  {}^p H^{-\dim(C)} Rg_* : \ Perv(-,\C) \to Perv(-,\C) $
is left exact.
Under $i^*$ and
$i_C^*$ perversity drops at most by $-d$. Notice $i_C^*\circ Rg_*
= Rg_*\circ i^*$ by proper base change. The distinguished  triangle
$({}^p H^{-\dim(C)}(Rg_*(K))[\dim(C)],Rg_*(K),K')$ for perverse $K$ defined by perverse truncation gives $K'\in {}^p D^{> -\dim(C)}(Z/C,\C)$ and
$i^*_C(K')\in {}^p D^{> -\dim(C)-d}(Z/C,\C)$. Therefore
$$ {}^p H^{-d}(i^*_C({}^p H^{-\dim(C)}(Rg_*(K)))\ \cong\ {}^p H^{-\dim(C)-d}(i^*_C Rg_*(K)) \ .$$
This gives both the upper (and lower) part of the 
commutative diagram
$$ \xymatrix@+0,8cm{ Perv(Z,\C) \ar[dr]^{{}^pH^{-\dim(C)-d}\circ i_C^*\circ Rg_*} \ar[r]^-{ {}^pH^{-\dim(C)}\circ Rg_*}\ar[d]_{{}^pH^{-d}\circ i^*}  &  \  Perv(Z/C,\C) \ar[d]^{{}^pH^{-d}\circ i_C^*} \cr  Perv(F,\C) \ar[r]^-{ {}^pH^{-\dim(C)}\circ Rg_*} & \ Perv(F/C,\C) }
$$  
so that
$$  {}^p H^{-\dim(C)}(Rg_* ({}^p H^{-d}i^*(K_\chi))) =  {}^p H^{-\dim(C)}(Rg_*({M^d})) \neq 0 $$
implies  ${}^p H^{-d}(i_C^*({}^p H^{-\dim(C)}Rg_*(K_\chi)))\neq 0$, and hence as required for step 2)
$$ {}^p H^{-\dim(C)}(Rg_*(K_\chi)) \neq 0 \ .$$

\medskip
Step 4) For the proof of the proposition for generic characters $\chi$ now suppose $$Rf_*(K_\chi)_b=0\ \quad  \quad \mbox{ at } \ b\in Z/A \ .$$     
Since $ M^{d} \neq 0 $, step 3 in the case $C=A$ therefore
implies that ${M}$
is acyclic on $F$, but not zero.
For a suitable $\chi$ all perverse cohomology sheaves $ M^i = {}^pH^{i}({M}) $ are acyclic on $F$ and their perverse Jordan-H\"older constituents as well.
Fix one such $\chi$ and notice ${}^p H^{-i}({M}_{\chi_0}) = {}^p H^{-i}({M})_{\chi_0} = M^i_{\chi_0}$ for any character $\chi_0$. 

\medskip
Step 5)  Since $F(A)=\emptyset$ by assumption, all the acyclic perverse Jordan-H\"older constituents of the perverse sheaves $M^i$  are negligible perverse sheaves on  $F\cong A$. Let $S\hookrightarrow  M^d \neq 0$ be a simple perverse subobject. Then $S=(\delta_C)_{\chi_0^{-1}} * L$ 
for some $L\in M(A)$  and some nontrivial abelian subvariety $C$ of $A$ and some character $\chi_0$. For the corresponding projection $g:X \to X/C$
$$  {}^p H^{-\dim(C)}(Rg_*(S_{\chi_0}))\ \ \hookrightarrow  
\ \ {}^p H^{-\dim(C)}(Rg_*(M^d_{\chi_0})) \ $$
by the left exactness of ${}^p H^{-\dim(C)}$. Furthermore ${}^p H^{-\dim(C)}(Rg_*(S_{\chi_0})) \neq 0 $.
For the last assertion notice
$Rg_*(S_{\chi_0}) = H^\bullet(C) \otimes_{\C} 
Rg_*(L)$ and hence $${}^p H^{-\dim(C)}( Rg_*(S_{\chi_0})) =
H^{-\dim(C)}(C) \otimes_\C \bigl( {}^p H^0(Rg_*(L)) \bigr)\ .$$  Indeed $L$ has nonvanishing Euler characteristic and therefore
$Rg_*(L) \neq 0$.  Then  ${}^p H^i( Rg_*(L))=0$ holds for all $i\neq 0$, since ${}^p H^i(Rg_*(S_{\chi_0}))=0$ for $i<-\dim(C)$. This being said, we obtain
$$ {}^p H^{-\dim(C)}(Rg_*(M^d_{\chi_0})) \neq 0 \ .$$
The proposition follows from step 1)-3) applied
for $\chi\chi_0$ instead of $\chi$. \qed

\medskip
{\bf Remark}.  For an abelian variety $X$, let $A$ be the connected stabilizer of an irreducible subvariety $Z$ of $X$, and let $\tilde Z$ be  the image of $Z$ in $\tilde X=X/A$. Then the connected stabilizer
of $\tilde Z$ in $\tilde X$ is trivial.

\begin{Proposition} \label{alsoincharp}
If $F(B)=\emptyset$ holds for all quotient abelian varieties  $B\neq X$ of $X$, then $$ F(X)=F_{max}(X)\ .$$ 
\end{Proposition}

{\it Proof}. Suppose $K\in F(X)$, but $K\notin F_{max}(X)$. Then there exists a minimal
quotient $p:X\to B$ such that $Rp_*(K_\chi)=0$ holds ($\chi$ generic). The fiber perverse
sheaves $M^i(b), b\in B$ and their Jordan-H\"older constituents $P$ are  
acyclic perverse sheaves on the abelian variety $Kern(p)$. Hence by the induction assumption these $P$ are in $N(Kern(p))$,
hence have degenerate support. Since $Z=supp(K)=\bigcup_{b\in B,i} supp(M^i(b))$, by 
[A] therefore $Z$ is degenerate, i.e. there exists an abelian subvariety $A$ of $X$ of
dimension $>0$ such that $Z+A=Z$. Suppose $A\neq X$.  For the quotient morphism $f:X\to X/A$
and by the induction assumption $F(A)=\emptyset$, then $Rf_*(K_\chi)\neq 0$ holds for generic
$\chi$ since otherwise $K\notin F(X)$ (proposition \ref{fantom}). Since  
$Rf_*(K_\chi)$ is perverse for generic $\chi$, therefore $L=Rf_*(K_\chi) \in E(X/A)$.
Furthermore $supp(L)=\tilde Z$ for $\tilde Z = f(Z)$ by proposition \ref{fantom}; in particular
$\tilde Z$ is irreducible. Since $L\in E(X/A)$ and since $F(X/A)=\emptyset$ holds by the induction assumption, the support $\tilde Z$ of $L$ is a finite union of degenerate subvarieties.    
Since $\tilde Z$ is irreducible, therefore $\tilde Z$ is degenerate. This is a contradiction because by  the remark above, for the quotient $f:X\to B=X/A$ by the connected stabilizer $A$,  the irreducible variety $\tilde Z=f(Z)$ has trivial connected stabilizer. 
This shows $Z=X$ and $Z$ is invariant under $Kern(p)$, so the vanishing $Rp_*(K_\chi)=0$
for generic $\chi$ implies $K\notin F(X)$ by proposition \ref{fantom} applied to 
the morphism $p: X\to B$.  \qed

\begin{Corollary} \label{2}  
If  $X$ is isogenous to the product $A_1\times A_2$
of simple abelian varieties $A_1,A_2$  of dimension $\dim(A_1)=\dim(A_2)=d$ (and $F(A_1), F(A_2)=\emptyset$), then $$ F(X)=F_{max}(X)\ .$$ 
\end{Corollary}

{\bf Remark}. In the situation of the last corollary
for irreducible $K\in E(X)$
either $\nu_K=d$, or $K\cong \delta_{X}^\psi$. 
Indeed $K\in F_{max}(X)$. So  
$\nu_K=d$ by lemma \ref{Kill1}. For $K\in N(X)$ we know $\nu_K=d$ or $2d$.

\medskip
Although not needed for the rest of the paper, 
in the remaining part of this section we show that under similar conditions
for $K\in F(X)$ the support of $K$
is $X$. We first state a  result of [W]. 

\begin{Theorem} \label{Support} For complex abelian varieties
the support $Z$ of an irreducible  perverse sheaves $K\in E(X)$ always is a degenerate irreducible subvariety of $X$, i.e. $A+Z=Z$ holds for some abelian subvariety $A$ of $X$ of dimension $>0$. 
\end{Theorem}

\begin{Corollary}\label{abram}  $K\in E(X) \Longrightarrow K\in N(X)$ for irreducible  $K\in Perv(X,\C)$, provided $F(A)=\emptyset$
holds where $A$ is the connected stabilizer of the support $Z$ of $K$.  
\end{Corollary}

{\it Proof}. 
By theorem \ref{Support} the support $Z$ of $K$ is an irreducible degenerate subvariety of $X$.  For the quotient $f:X\to B=X/A$ by the connected stabilizer $A$,  the irreducible variety $\tilde Z=f(Z)$ by construction has trivial connected stabilizer. If for generic $\chi$ the perverse sheaf $Rf_*(K_\chi)$ vanishes,  our assertion  for $K\in F(X)$ and $F(A)=\emptyset$ is a consequence of  proposition \ref{fantom}.  So it suffices to show $Rf_*(K_\chi)=0$. Suppose
$Rf_*(K_\chi)\neq 0$. Then  $\tilde Z= supp(Rf_*(K_\chi))$ by proposition \ref{fantom}. 
Since $\chi(Rf_*(K_\chi))=\chi(K)=0$, any irreducible constituent 
$P$ of $Rf_*(K_\chi)$ is in $E(X/A)$ and therefore $supp(P)$ is degenerate by 
theorem \ref{Support}. 
Since $\tilde Z=supp(Rf_*(K_\chi))=\bigcup supp(P)$ and $\tilde Z$ is irreducible, we obtain $\tilde Z = supp(P)$ for some $P$. Hence $\tilde Z$ is degenerate. A contradiction. \qed

\medskip
An immediate consequence of these arguments is

\begin{Proposition}\label{Rec} Assume $F(B)=\emptyset$ for all quotients $B$ of $X$ of dimension $< \dim(X)$. Then any $K\in F(X)$ 
has support $X$. 
\end{Proposition}

{\it Proof}. The assertion $supp(K)=X$ follows from corollary \ref{abram}. To show $K\in F_{max}(X)$ consider  $f:X\to B$ for minimal $B$.
Suppose $Rf_*(K_\chi)=0$ holds for generic $\chi$. Then proposition \ref{fantom} can be applied for $K$ and $A=Kern(f)$, since $F(A)=\emptyset$
and $supp(K)=X$ is invariant under $A$. This proves $K\in N(X)$.
A contradiction. \qed

\medskip

\medskip

\section{Proof of the main theorem} \label{main}

\medskip

\begin{Lemma}\label{B} $F_{max}(X)=\emptyset$ if
 $X$ has a simple quotient  $B$ 
with $\dim(B)>\mu(X)$.   
\end{Lemma} 

\medskip
{\it Proof}. Suppose $K\in F_{max}(X)$. 
For  $f:X \to B$ then $\mu(M^i)\geq \mu(X)$
for $i=0,..,\dim(B)-1$ by lemma \ref{A}, since these $M^i$ are acyclic. 
Therefore 
$ {\cal H}^{-\mu(X)} K_\chi\vert_{f^{-1}(0)} \cong 
{\cal H}^{-\mu(X)}({M}) \cong {\cal H}^{-\mu(X)}(M^0) = 0 $. Indeed,
for $K\in F_{max}(X)$ the support  is not 
contained in a translate of a proper 
abelian subvariety of $X$ by lemma \ref{maximal-}, so $M^0=0$. 
This shows
$$  \mu(K) > \mu(X) \ .$$ Now also ${\cal P}_K\in F_{max}(X)$ by lemma \ref{Kill1}. The last inequality applied for  ${\cal P}_K$ instead of $K$ leads to a contradiction. Indeed  
 $\mu({\cal P}_K)= \nu_K = \mu(X)$ holds by the maximality of $K$  using lemma \ref{Kill1}. \qed

\medskip
{\it At least three simple constituents}.
Consider quotients $g:X \to B$ with nontrivial  
kernel $C$, where $B$ has two simple factors $B_1$ and $B_2$ and
$$ 0 \to B_1 \to B \to B_2 \to 0 \ .$$ 
By lemma \ref{B} we may assume that all simple
factors of $X$ have the same dimension $d=\mu(X)$. Put
$A= g^{-1}(B_1)$.  
Consider $K\in F(X)$. Then
for $f:X\to X/A=B_2$ and
generic $\chi$ 
the perverse sheaf $Rg_*(K_\chi)$ is in $E(B)$, hence bx the decomposition theorem of  the form  
$$    Rg_*(K_\chi) \ = \ \bigoplus_{i\in I} \  T^*_{b_i}(\delta_{A_i}^{\psi_i})*M_i 
\quad \oplus \ \ \ \mbox {objects in } F(B) \ \ \oplus \quad rest  \ $$
for finitely many  simple abelian subvarieties $A_i$  of $B=B_1\times B_2$ 
of dimension $ \dim(A_i)= \mu(X) = d$. Here
$b_i$ are certain points in $B$, $M_i$ are in  $M(B)$ and $\psi_i$ are certain 
characters with $\delta_A^\psi := \C_A[\dim(A)]_\psi$. The term \lq{rest}\rq\ denotes the $B$-invariant term. 

\medskip
{\bf Claim}. {\it For $K\in F_{max}(X)$  the index set $I$ is empty.
}

\medskip
{\it Proof}. Assume $I$ is not empty. For generic $\chi$
the perverse sheaf $Rg_*({\cal P}_{K,\chi})$ admits for each $i\in I$ a 
nontrivial morphism  in the derived category ([W2, prop.1])
$$  Rg_*({\cal P}_{K\chi})[\nu_K] \to {\cal P}_{P_i}[\nu_{P_i}] $$
for each summand $P_i =   T^*_{b_i}(\delta_{A_i}^{\psi_i})*M_i$.
Since $\nu_K=\mu(X)=d$ (by lemma \ref{Kill1} this holds for maximal $K$
)
and also $\nu_{P_i}=d$
for $i\in I$, these morphisms define nontrivial morphisms of perverse sheaves
$$ Rg_*({\cal P}_{K\chi})\ \longrightarrow\ {\cal P}_{P_i} =  \delta_{A_i}^{\psi_i} \ ,$$
which then are epimorphisms. 
Therefore by the decomposition theorem
$$ Rg_*({\cal P}_{K,\chi}) \ = \ \bigoplus_{i\in I} \  \delta_{A_i}^{\psi_i} \ \ \ \oplus \ \ \ others \ .$$
Since $P={\cal P}_{K,\chi}$ also satisfies $\nu_P=\nu_K=d$ and also is in $F(X)$, 
we may replace $K$ by ${\cal P}_K$. Thus it suffices to
show that no terms $ \delta_{A_i}^{\psi_i}$ can appear in $Rg_*(K_\chi)$
under our assumptions on $K$ above.
Indeed, 
proposition \ref{gut} implies $$R^{-d}g_*(K_\chi)_{b'_1,b'_2}=0$$ 
for almost all $b'_1\in B'_1$ (for fixed $b'_2\in B'_2$). 
For this we have to choose the map $f$ used in this proposition \ref{gut}, so that
$B_1$ from that proposition is $B'_1=g^{-1}(A_i)$ for some fixed $i$ (and not our 
fixed $B_1$). Then for $b'_2=0$ the fiber $$R^{-d}g_*(K_\chi)_{b'_1,0}={\cal H}^{-d}( \delta_{A_i}^{\psi_i})_{b'_1,0} $$ is not zero for almost all $b'_1\in A_i$ contradicting proposition \ref{gut}, if $I\neq \emptyset$. This proves the claim. In fact,
the support and vanishing 
assumptions from page \pageref{ass} are now satisfied for ${\cal P}_K$ with $d=\mu(X)$.
Recall that the support  ${\cal P}_K$ must not lie in a proper abelian subvariety for 
 proposition \ref{gut}.  But $K\in F_{max}(X)$ implies ${\cal P}_K\in F_{max}(X)$  
by lemma \ref{Kill1}; hence lemma \ref{maximal-} takes care of this. \qed

\medskip
From the preceding dicussion we conclude 

\begin{Lemma} \label{strange} Suppose given $g:X\to X/A=B$ with $B$ and $A\neq 0$ as above. 
Then for $K\in F_{max}(X)$ and generic $\chi$ we have
$$ \fbox{$  Rg_*(K_\chi) = \mbox {objects in } F(B) \ \oplus \ \bigoplus_{\psi \in \Psi(\chi)} \ m_\psi \cdot \delta_{B}^{\psi} $} \ .$$
\end{Lemma}

\begin{Theorem} \label{MAINTHEOREM} $F(X)=\emptyset$ holds for all complex abelian varieties $X$.
\end{Theorem}

{\it Proof}. We show $F(X)=\emptyset$  by induction on the number of simple factors (or the dimension) of $X$.  Recall, for the simple case this was shown in [KrW]. The cases with two simple factors will 
be considered in proposition \ref{ToDo} below. So assume that $X$ has at least three simple 
factors and that $F(D)=\emptyset$ already holds for all proper subvarieties or proper quotients $D$ of $X$.   
Hence 
$K\in  F_{max}(X)$ by proposition \ref{alsoincharp}. 
Then by lemma \ref{B} all simple factors have dimension $\mu(X)=d$, so there  exists a quotient $g:X\to B$ with kernel $A\neq 0$ 
$$ 0 \to B_1 \to B \to B_2 \to 0 \ $$ 
with simple factors $B_1$ and $B_2$ of dimension $d=\mu(X)$.
Now $F(B)=\emptyset$ holds by the induction assumption and $K$ is maximal.
 Hence  $Rg_*(K_\chi)= \bigoplus_{\psi \in \Psi(\chi)} \ m_\psi \cdot \delta_{B}^{\psi}$ is $B$-invariant by lemma \ref{strange} for generic $\chi$. But this contradicts
 the maximality of $K$, since
then $Rf_*(K_\chi)=0$ vanishes for generic $\chi$. \qed  

%

\medskip

\section{Functoriality}\label{previous}


\begin{Lemma} \label{projection} 
For simple $A=  Kern( p: X \to B)  $ and for $K\in F(X)$ the perverse cohomology sheaves ${}^p H^i(Rp_*(K))$ are in $ E(B)$ for all $i$.
\end{Lemma}

{\it Proof}.  Using isogenies one reduces this to the case $X=A\times B$, where  $p$ is the projection onto the second factor. Then either $Rp_*(K)= 0$ and there is nothing to prove, or for certain irreducible perverses sheaves $P_i$ and certain integers
$\nu_i $ 
$$   Rp_*(K) = \bigoplus_{i} \  P_i[\nu_i] \quad , \quad \vert \nu_i \vert < \dim(A)\ .$$
Suppose for one of the irreducible summands $P_i\notin E(B)$. For the projection $q$ onto the first factor $A$
$$ \xymatrix@+0,5cm{ \ X=A\times B \ \ar[r]^-{p} \ar[d]_q^-\chi &   \ B\  \ar[d]^-{\chi} \cr
\ A\ \ar[r] & \ Spec(\C) \ } $$
by lemma \ref{relGVT}  there exists a (generic) character $\chi$ of $\pi_1(X,0)$ 
so that both 1) $Rq_*(K_\chi)$ is perverse and 2) $H^\bullet(B,(P_i)_\chi)=H^0(B,(P_i)_\chi)$
holds for the finitely many $P_i$.
Since $K\in F(X)$, the Euler characteristic of $Rq_*(K_\chi)$ vanishes. Since $Rq_*(K_\chi)$ is perverse, therefore $Rq_*(K_\chi)\in E(A)$ is 
either zero or of the form 
$ Rq_*(K_\chi) = \bigoplus_{\psi\in \Psi} \ m_\psi\cdot \delta_A^{\psi} $
and
$   H^\bullet(A,Rq_*(K_\chi)) \ =\
 \bigoplus_{i} \ H^\bullet(B,(P_i[\nu_i]_\chi)$. 
For $P_i\notin E(B)$ the cohomology of all twists $(P_i)_\chi$ does not vanish
(since the Euler characteristic is constant $>0$ and independent of twists). So for generic $\chi$ then $   H^\bullet(B,P_i[\nu_i]_\chi) \ = \ H^{\nu_i}(B,(P_i)_\chi) \neq 0 \ .$
Hence by comparison $Rq_*(K_\chi)$ can not vanish
and therefore is a sum of $\delta_A^{\psi}$ so that
at least for one character $\psi$ the cohomology $H^\bullet(A,\delta_A^{\psi})$
does not vanish, i.e. $\psi$ is trivial.
This gives a contradiction, since then $H^\bullet(A,\delta_A^{\psi})$
contains terms of degree $\dim(A)$. For all summands $P_j$, that are in $E(B)$, the cohomology $
 H^\bullet(B,P_j[\nu_j]_\chi)$ vanishes for generic $\chi$ by corollary \ref{azyk}.  For the others the 
 cohomology $H^{\nu_i}(B,(P_i)_\chi) \neq 0 $ 
does not contribute to degree $\dim(A)$, since these $\nu_i$ satisfy 
$\vert \nu_i\vert < \dim(A)$.  \qed

\section{Convolution}\label{next}

For simple  abelian varieties $A_1$ and $A_2$ of dimension $d$ and simple
perverse sheaves $K,L$ on the cartesian product $X=A_1\times A_2$ consider the diagram
$$ \xymatrix@+0,3cm{ A_{1}\times A_{1} \ar[d]^a  &   X\times X   \ar[d]^b  \ar[l]_{p_1 \times p_1}\ar[r]^{p_2 \times p_2} &   A_{2} \times A_{2}\ar@{=}[d] \cr
A_{1}   &   Y=A_{1}  \times A^2_{2}  \ar[d]^c \ar[l]_{p_1 }\ar[r]^{p_{23}} &   A_{2} \times A_{2} \ar[d]^a\cr  A_{1} \ar@{=}[u] 
& A_{1}\times A_{2} \ar[r]^{p_2} \ar[l]_{p_1} & A_{2}  }$$
with the morphisms
$a(x,y)=x+y$,
$ b(x_1,x_2,x_3,x_4)=(x_1+x_3,x_2,x_4)$ and
 $c(y_1,y_2,y_3) = (y_1,y_2+y_3) $
and the projections $ p_{23}(y_1,y_2,y_3) = (y_2,y_3) $ resp.
$  (p_1\times p_1)(x_1,x_2,x_3,x_4) = (x_1,x_3) $ and
$  (p_2\times p_2)(x_1,x_2,x_3,x_4) = (x_2,x_4) $ and $Y:= A_1 \times A_2^2$.
Then $$ K{\circledast} L := Rb_*(K\boxtimes L) \quad , \quad
  Rc_*(K\circledast L) = K*L \ .$$
By the decomposition theorem  $K\circledast L$ is a semisimple complex.
By the relative K\"unneth formula $   R(p_i\times p_i)_*(K\boxtimes L) =  Rp_{i*}K \boxtimes Rp_{i*}(L)$ and hence 
$$ Rp_{23*}(K\circledast L) =  Rp_{2*}(K) \boxtimes Rp_{2*}(L) \ .$$

\medskip
By twisting both perverse sheaves $K,L$ with the same character 
$\chi=(\chi_1,\chi_2)$ of $\pi_1(X,0)=\pi_1(A_1,0)\times \pi_1(A_2,0)$ 
the direct images $Rp_{2*}(K)$ and $Rp_{2*}(L)$, and also
$Rp_{23*}({}^p H^i(K\circledast L))$ become
 perverse sheaves on $A_2^2$ for generic $\chi_1$ (lemma \ref{relGVT}).   
Now $$  \bigoplus_{i+j=k} {}^pH^j(Rp_{23*}({}^pH^i(K\circledast L)))  =  {}^pH^k(Rp_{23*}(K\circledast L)) $$ holds by the decomposition theorem.
The right side vanishes for $k\neq 0$, since $ Rp_{23*}(K\circledast L) =  Rp_{2*}K \boxtimes Rp_{2*}(L) $ is perverse. So all terms for $k\neq 0$ on the left are zero, the terms for $j\neq 0$ vanish after twisting with a suitable character $\chi_1$. Hence
$ Rp_{23*}({}^pH^0(K\circledast L))) =  Rp_{23*}(K\circledast L)  =  Rp_{2*}(K) \boxtimes Rp_{2*}(L) $ and
$ Rp_{23*}({}^pH^i(K\circledast L))) =  0$ for  $i \neq 0  $ and
generic $\chi_1$.
Notation: $K\circ L = {}^pH^0(K\circledast L)$ and $K\bullet L = \bigoplus_{i\neq 0} {}^pH^i(K\circledast L)$. Since 
$\ \   K\circledast L = (K\circ L) \ \bigoplus \ (K\bullet L)\ \ $, we obtain
$$   Rp_{23*}(K\circ L) = Rp_{2*}(K)\boxtimes Rp_{2*}(L) \quad , \quad 
     Rp_{23*}(K\bullet L) = 0 \ .$$
Notice $Rp_{23*}(P)=0$ for simple constituents $P$ of the semisimple perverse sheaf 
$K\bullet L$,  hence 
$P\in E(Y)$. 
Since $Kern(p_{23})=A_1$ is simple and $Rp_{23*}(P)=0$, lemma \ref{notzero}
implies  $P\in N(Y)$. Then there exists an abelian subvariety $B$ and a  character $\psi$ and some $M\in M(Y)$ such that $P=\delta_B^\psi*M$ holds, and $Rp_{23*}(P)
=Rp_{23*}(\delta_B^\psi)*Rp_{23*}(M)=0$ implies $Rp_{23*}(\delta_B^\psi)=0$
so that $B \subset Kern(p_{23})$. Hence  up to a character twist $P=p_{23}^*[d](Q)$ for some 
$Q\in Perv(A_2^2,\C)$. 
Indeed $Stab(P)^0$ is an abelian subvariety of $Y$, and therefore
$Rp_{23*}(P_\chi)=0$ for generic $\chi$ implies $Kern(p_{23})\subseteq Stab(P)^0$.
Hence

\begin{Lemma} \label{technical} Up to character twists the irreducible constituents 
of $K\bullet L$ are in $ p_{23}^*[d](D_c^b(A_2^2,\C))$.
\end{Lemma}

\bigskip\noindent
By proper base change we obtain

\begin{Corollary} \label{hilfe} Up to character twists the irreducible constituents 
of 
$Rc_*(K\bullet L)$ are in $ p_{2}^*[d](D_c^b(A_2,\C))$.
\end{Corollary}

\bigskip\noindent
By cohomology bounds for the perverse cohomology of $K*L$
we obtain
  
\begin{Corollary}
For $K,L\in Perv(X,\C)$ up to character twists the irreducible constituents of
$ {}^p \tau^{>d}(K*L)$ are in $p_2^*[d](D_c^b(A_2,\C))$.
\end{Corollary}

{\it Proof}. For the proof we may twist both $K$ and $L$ by an arbitrary character
$\chi=(\chi_1,\chi_2)$ of $\pi_1(X,0)$, since convolution and therefore also $b$ and $c$ commute with 
character twists. Hence the claim follows
from $Rc_*(K\circ L) \subseteq {}^pD^{[-d,d]}(X)$ and corollary \ref{hilfe}. \qed

\medskip
The role of the indices $1$ and $2$ of the decomposition $X = A_1\times A_2$ is arbitrary, so by a switch

\begin{Corollary}
For $K,L\in Perv(X,\C)$ and $\vert i\vert >d$  
 we have $ {}^pH^i(K*L)  =  
\bigoplus_\varphi\ m_{i\varphi}\cdot \delta_{X}^\varphi $.
\end{Corollary}

\begin{Corollary} \label{funny}
For $K,L\in F(X)$ all the summands $P[d] \hookrightarrow K*L$ for which $P$ is not $X$-invariant appear in the form $H^\bullet(\delta_{A_2})\otimes_\C P \hookrightarrow K*L$ so that for some $P\in Perv(X,\C)$ and $Q\in {}^p D^{[-d+1,d-1]}(X)$ the following holds
$$ K*L \ = \ \Bigl(H^*(\delta_A)\cdot P\Bigr)\ \oplus\ \
\bigoplus_{i,\varphi \in \Phi(\chi)} \ m_{i\varphi } \cdot \delta_{X}^\varphi[-i] \
 \ \oplus\ \ Q  \ .$$
\end{Corollary}
 
{\it Proof}. First assume only  $K,L\in Perv(X,\C)$. Then any term $P[d]$ in $K*L$ with perverse $P$, not invariant under $Kern(p_2)$, is in $Rc_*(K\circ L)$; since $K\circ L\in Perv(A_1\times A_2^2,\C)$
and since $c$ is smooth of relative dimension $d$, hence $c^*(P)[d]$ must be a summand of $K\circ L$ and then
$Rc_*(c^*(P)[d])=H^\bullet(\delta_{A_2})\cdot  P$ is a summand of $Rc_*(K\circ L)$. 
The same applies for terms $P[d]$ in $Rc_*(K\circ L)$ that are 
$Kern(p_2)$-invariant. But now there might be $Kern(p_2)$-invariant
terms $P[d]$ not coming from  $Rc_*(K\circ L)$ but only from $Rc_*(K\bullet L)$.
To exclude this possibility we assume $K,L\in F(X)$ so that we can apply the next
corollary to show that these critical summands are $X$-invariant, hence contained in $\bigoplus_{i,\varphi \in \Phi(\chi)} \ m_{i\varphi } \cdot \delta_{X}^\varphi[-i]$.
\qed 

\medskip
We now allow arbitrary characters $\chi_1$. Notice 
$Rb_*(K_{\chi_1}\boxtimes L_{\chi_1}) = Rb_*(K\boxtimes L)_{\chi_1}$ or
$(K\circledast L)_{\chi_1} = (K_{\chi_1}\circledast L_{\chi_1})$. 
For irreducible $P=p_{23}^*[d](Q)$ we can determine $Q$ as a direct summand of $Rp_{23*}(P_{\chi_1})$ for a suitable choice of
$\chi_1$. Now $Rp_{23*}(P_{\chi_1}) = Rp_{2*}(K_{\chi_1})\boxtimes 
Rp_{2*}(L_{\chi_1})$.
For $K,L\in F(X)$ and arbitrary $\chi_1$ we then get $Rp_{2*}(K_{\chi_1}) = \bigoplus_{i,\psi}
m_{\psi i}\cdot \delta_{A_2}^\psi[-i]$ from lemma \ref{projection}, since $A_1$ is simple.
Then
$Q$ has to be $A_2^2$-invariant for $K,L\in F(X)$. This proves 

\begin{Corollary}\label{volli}
For $K,L\in F(X)$ the complex $K\bullet L$ is $Y$-invariant.
\end{Corollary}

\medskip
Similarly $Rp_{1*}(K)$ and  $Rp_{1*}(L)$ are perverse,  by twisting with 
a generic $\chi_2$   
and 
$ \bigoplus_{i+j=k} {}^pH^j( Rp_{1*}({}^p H^i(K\circledast L))) \ = \ {}^p H^k(Rp_{1*}(K)*Rp_{1*}(L))$. 
Indeed from now on we make the 

\medskip
{\bf Assumption}. $K,L\in F(X)$.

\medskip
Then as required: 1) For generic $\chi$ we have $Rp_{2*}(K)$ and 
$Rp_{2*}(L)$ are perverse and translation invariant on $A_2$; and similarly  
$Rp_{1*}(K)$ and 
$Rp_{1*}(L)$ are perverse and translation invariant on $A_1$.
2) Furthermore ${}^pH^i(K\circledast L)=0$ for $\vert i\vert \geq d$.
Of course only the case $i=\pm d$ is relevant. If 
${}^pH^{\pm d}(K\circledast L)\neq 0$, then $K\boxtimes L$ is $Kern(b)$ invariant. Hence $K\boxtimes L$ does not lie in $F(X^2)$, contradicting the fact that $K,L\in F(X)$ implies $K\boxtimes L \in F(X^2)$.
Indeed, suppose $A\hookrightarrow X^2$ is a nontrivial abelian subvariety and $K\boxtimes L$ would be invariant under $A$. Then $K$ is invariant under $pr_1(A)$ and $L$ is invariant under $pr_2(A)$. 
By our assumptions, hence $pr_1(A)=0$ and $pr_2(A)=0$. A contradiction.


\medskip\noindent
By  lemma \ref{projection} for generic $\chi_2$ 
(we do not write this twist !)
$$  Rp_{1*}(K) = \bigoplus_\psi \ m_\psi^K \cdot \delta_{A_1}^\psi \quad , \quad
  Rp_{1*}(L) = \bigoplus_\psi\  m_\psi^L \cdot \delta_{A_1}^\psi \ $$
and hence
$$  Rp_{1*}(K\circledast L) = H^\bullet(\delta_A) \cdot
\bigoplus_\psi \ m_\psi^Km_\psi^L \cdot \delta_{A_1}^\psi \ .$$
Since $Rp_{1*}(K\circledast L)= Rp_{1*}\circ Rc_*(K\circledast L)$, we now compare this 
with  $Rp_{1*}(K*L)= Rp_{1*}\circ Rc_*(K\circledast L)$
using the formula for $K*L$ obtained in corollary \ref{funny} 
$$Rp_{1*}(K*L) \ =\  Rp_{1*}\Bigl( (H^*(\delta_A)\cdot P)\ \oplus 
\bigoplus_{i, \varphi\in \Phi(\chi)} \ m_{i\varphi} \cdot \delta_{X}^\varphi[-i] \ 
\oplus\  Q\Bigr) $$ $$ =H^\bullet(\delta_A)\cdot  \Bigl(
Rp_{1*}(P) \oplus \bigoplus_{\varphi'} \ m_{0\varphi'} \cdot \delta_{A_1}^{\varphi'} \Bigr)\ \oplus\ \
Rp_{1*}  (Q) \ .$$
The terms $Rp_{1*}(\delta_{X}^\varphi[-i] )$ must vanish for $i\neq 0$, 
since they have no counter part in the comparison! 
But for $i=0$ there might be contributions $Rp_{1*}(\delta_{X}^\varphi[-i]) =
\delta_{A_1}^{\varphi'}$ for certain $\varphi$.  
Notice, for generic $\chi_2$ the functor $Rp_{1*}$ preserves 
perversity. Therefore, making a comparison of the terms of degree $d$ in $H^\bullet(\delta_A)$ (for $A=A_1$ or $A_2$)
first, we immediately get $Rp_{1*}(Q)=0$ for generic $\chi_2$.
Hence for the constituents of $Q$ are up to twists in $p_1^*[d](Perv(A_1,\C))$.
For those in $N(X)$ the assertion is obvious. For those in $F(X)$ use corollary \ref{2}, $F(A_2)=\emptyset$ and
proposition \ref{fantom}.
Since the decomposition $X= A_1\times A_2$ is arbitrary, 
again by switching the  indices  we obtain the next

\begin{Proposition}\label{deco}
For $K,L\in F(X)$ there exists $P\in Perv(X,\C)$
and some $X$-invariant semisimple complex $T\in D_c^b(X,\C)$ so that
$$ \fbox{$  K*L \ = \ H^*(\delta_{A_2}) \cdot P \ \oplus \ T $} \quad \quad \fbox{$ H^\bullet(\delta_{A_2}) \cdot P \ \hookrightarrow\
Rc_*(K\circ L) $}  \ ,$$ 
and $H^\bullet(\delta_{A_2}) \cdot P = Rc_*(c^*(P)[d]) \subseteq Rc_*(K\circ L)$. 
If $K,L$ are monoids, then we have $P=K$ for $K\cong L$ respectively $P=0$
for $K\not\cong L$.
\end{Proposition}

\medskip
Concerning the notation used: We call a monoidal irreducible perverse sheaf on $X$ that is contained in $E(X)$  a monoid.

\medskip
{\it Proof}. 
From our discussion it is clear that $P=P_1 \oplus P_2$ decomposes
such that $P_1\in Rc_*(K\circ L)$ and $P_2$ by corollary \ref{volli}
 is $X$-invariant. This allows to replace $P$ by $P_1$ and $T$ by $T\oplus H^\bullet(\delta_{A_2})\cdot P_2$. The assertion on the precise form of
 $P$ for monoids $K$ and $L$ (see [W2]) then follows from [W2, lemma 7] applied to the class $F(X)$
 after localization with respect to the hereditary class of $X$-invariant complexes. \qed

\section{Sheaves on $X=A\times B$}  \label{specialsheaves}

\medskip
For a monoidal perverse sheaf $K\in F(X)$ on $X=A\times B$, for simple abelian varieties $A,B$ with dimension $\dim(A)=\dim(B)$, in the last section we have shown that $K\cong K^\vee$ and $K*K \cong H^\bullet(\delta_B)\cdot K\oplus T$ holds for some
$T= \bigoplus_{\varphi} m_{i\varphi} \delta_X^\varphi[-i]$.  Let $p:X \to B$ denote the projection onto the second factor.

\medskip
Replacing $K$ by a twist $K_\chi$ we can suppose  
that $m_{i\varphi}=0$ holds for the trivial character $\varphi$. We can furthermore assume that there exists  an abelian subvariety $A\hookrightarrow X$ with quotient $p:X\to B=X/A$ such that all the finitely many $\varphi$ with $m_{i\varphi}\neq 0$ have nontrivial restriction on $A$. Assuming this, then $Rp_*(T)=0$ and therefore
$Rp_*(K)^2 = H^\bullet(\delta_B) \cdot Rp_*(K)$. Then $L:=Rp_*(K) = \bigoplus m^K_{i\psi} \delta_B^\psi[i] $ by lemma \ref{projection}. Hence $Rp_*(K)$ must be  perverse ($m_{i\psi}=0$ for $i\neq 0$) so that $\bigoplus_\psi H^\bullet(\delta_B) \cdot m_\psi^2 \cdot \delta^\psi_B \cong
\bigoplus_\psi H^\bullet(\delta_B)\cdot m_\psi \cdot \delta^\psi_B$. Both statements are an easy consequence 
of $L*L \cong H^\bullet(\delta_B)\cdot L$. Therefore 
$$  Rp_*(K) \ = \ \bigoplus_{\psi} \ m_\psi\cdot \delta_B^\psi \quad , \quad m_\psi\in \{ 0,1\} \ .$$
Notice, that by replacing $K$ by some other $K_\chi$ with the same properties
the coefficients $m_\psi$ may of course change. But 
the generic rank $r$ $$r= \sum_{\psi} m_\psi =  rank(R^{-d}p_*(K))_b = (-1)^{\dim(A)}\chi(H^\bullet(A,M(b)))$$ is independent from the specific character twist and the point $b\in B$, since
$$\chi(H^\bullet(A,M(b)))=\chi(H^\bullet(A,M(b)_{\chi_0}))$$ holds for every (!) character $\chi_0$.  Recall that $H^\bullet(A,M(b)_{\chi_0}) = Rp_*(K_{\chi_0})_b$ holds for every $b$ and every $\chi_0$. 

\medskip
$K$ and any of its twists $K_\chi$ have the same 
stabilizer. Under our assumptions the stabilizer $H$ of $K$ is a finite
subgroup of $X$, say of order $n$.  Then, for the isogeny $\pi=n \cdot_X: X \to X$, we may replace $K$ by one of the irreducible components $P$ of $\pi_*(K)=\bigoplus_{\chi\in H^*}  P_\chi$, which is a direct sum of $\# H$ simple monoidal perverse sheaves $P_\chi$ with $\nu_P=\nu_K$ and each $P\in F(X)$ has trivial stabilizer
[W2, cor. 4]. So in the following we may always assume ${\cal H}^{-d}(K)=\delta_0$, i.e. that $T_x^*(K)\cong K$ implies $x=0$. 

\begin{Theorem} \label{cp}  For simple complex abelian varieties $A$ and $B$
of dimension $\dim(A)=\dim(B)$ let denote $X=A\times B$ and $p:X\to B$ the projection
onto the second factor. For a fixed translation invariant complex $T$ on $B$ consider the set ${\cal B}$ of isomorphism classes of perverse sheaves $K$ on $X$ such that
\begin{enumerate}
\item $K$ is an irreducible perverse sheaf in $Perv(X,\C)$.
\item $K^\vee \cong K$.
\item $K*K \cong H^\bullet(\delta_B) \cdot K\ \oplus \ T$.
\item The support $supp(K)$ is not contained in a simple abelian subvariety.
\item $Stab(K)$ is trivial, i.e. $T_x^*(K)\cong K$ implies $x=0$.
\item $Rp_*(K)$ is perverse of generic rank $r>1$.
\end{enumerate}
Then ${\cal B}$ is empty.
\end{Theorem}

{\it Proof}. For simplicity of exposition we may assume $supp(K)=X$ (by prop.\ref{Rec}) although this is not essential for the argument. We prove this theorem by reducing it to a corresponding statement for base fields of positive characteristic $p$, using the method of Drinfeld [D] and [BK], [G]. The conditions defining ${\cal B}$ are constructible conditions in the sense of [D, lemma 2.5] and [D, section 3].
If ${\cal B}$ were not empty, the argument of [D] therefore provides us with some other $K\in {\cal B}$, which now is a ${\frak o}$-adic perverse sheaf on $X$ for some finite extension ring ${\frak o}$ of an $l$-adic ring $\Z_l$ with prime element $\pi$ generating the maximal ideal of ${\frak o}$, such that furthermore 
we find a subring $R\subset \C$ finitely generated over $\Z$
so that $X$ and and the complex $K\otimes^{\rm L}  {\frak o}/\pi {\frak o}$ is defined over $Spec(R)$, i.e. the pair $(X,K\otimes^{\rm L}  {\frak o}/\pi {\frak o})$ descends to some $(X_R,K_R^{(\pi)})$, with the following properties: a) The structure morphism $X_R \to S=Spec(R)$ is universally locally acyclic with respect to  $K_R^{(\pi)}$, b) For every maximal ideal of $R$ with the finite residue field $\kappa$ and the corresponding strict Henselization $V$ for a geometric point over this maximal ideal and embeddings $R\subset V \subset \C$, such that the conditions of [D, section 4.9]\footnote{similar to [BBD, lemma 6.1.9] where $V$ instead is chosen non canonically as some strict Henselian valuation ring with center in the maximal ideal.} are satisfied for a suitable ${M_i}$
attached to $K_R^{(\pi)}$ and its convolution square, one has equivalences of categories $$D_{\{M_i\}}^b(X,{\frak o})\ \sim\  D_{\{M_i\}}^b(X_R \otimes_R V,{\frak o})\ \sim\ D_{\{M_i\}}^b(\overline X,{\frak o}) \ $$
so that the structure morphism $f: X_V \to Spec(V)$ is universally locally acyclic with respect to  $K_V^{(\pi)}$, and c) The reduction
$\overline K$ of $K$ is an irreducible perverse Weil sheaf on the special fiber $\overline X$ with structure morphism $\overline f: \overline X \to Spec(\overline\kappa)$ defined over some finite extension of $\kappa$, and finally d) there are similar equivalences of categories as in b) 
 $$D_{\{N_i\}}^b(B,{\frak o})\ \sim\  D_{\{N_i\}}^b(B_R \otimes_A V,{\frak o})\ \sim\ D_{\{N_i\}}^b(\overline B,{\frak o}) \ $$
for suitable ${N_i}$  on $B_R$ attached to the perverse  cohomology sheaves of the sheaf complex $Rp_*(K_R^{(\pi)})$ and its convolution square similar to [D, 6.2.3] with a commutative diagram ($\bf *$) similar as in [D, (6.1)]
$$  \xymatrix{  D_{\{M_i\}}(X,{\frak o}) \ar[d]_{Rp_*} \ar[r]^\sim & D_{\{M_i\}}(\overline X,{\frak o}) \ar[d]^{R\overline {p}_*}    \cr
 D_{\{N_j\}}(B,{\frak o}) \ar[r]^\sim & D_{\{N_j\}}(\overline B,{\frak o}) } $$
Here $D_{\{M_i\}}(X,{\frak o})$ as a full subcategory of $D(X,{\frak o})$ is defined in [D, 4.9] as the inverse 2-limit of subcategories 
$D_{\{M_i\}}(X,{\frak o}/\pi^r{\frak o}) \subset D_{prf}(X,{\frak o}/\pi^r{\frak o})$
so that $C \in  D_{\{M_i\}}(X,{\frak o})$ iff $C\otimes^{\rm L}_{\frak o} {\frak o}/\pi{\frak o}$  
is in the thick triangulated subcategory  $D^{\{M_i\}}(X,{\frak o} /\pi{\frak o})$ of $D(X,{\frak o} /\pi{\frak o})$ generated by the $M_i $ in a finite set $\{M_i\}$ of fixed complexes
$M_i \in D_c^b(X_R,{\frak o}/\pi {\frak o})$. Similarly for $\overline X$ respectively for finitely many $N_j$ on $B$ (or $\overline B$). 
We briefly remark that in loc. cit. these equivalences of derived categories above over $\C,V,\overline \kappa$ are first proved 
on the level of ${\frak o}/\pi^r{\frak o}$-coefficients. There one implicitly uses that for perfect complexes $K$
of ${\frak o}/\pi^r{\frak o}$-sheaves (or projective limits of perfect complexes of such complexes i.e. in the sense of [KW, p.96f]) the distinguished triangles
$(K\otimes_{\frak o}{\frak o}/\pi^{r-s}{\frak o}, K\otimes_{\frak o}{\frak o}/\pi^{r}{\frak o},
K\otimes_{\frak o}{\frak o}/\pi^{s}{\frak o})$ for $0\leq s \leq r$ show that
$K\otimes_{\frak o}{\frak o}/\pi^{r}{\frak o}$ is contained in $D^{\{M_i\}}(X,{\frak o} /\pi^r{\frak o})$,
if $K\otimes_{\frak o}{\frak o}/\pi{\frak o}$ is contained in $\{ M_i\}$ or $D^{\{M_i\}}(X,{\frak o}/\pi{\frak o})$.

\medskip
The first condition a) can be achieved by a suitable localization of $R$ using
[Fin, thm. 2.13]. The acyclicity conditions in loc. cit. $i^*{\cal F} \cong i^*R\overline j_*{\cal F}$ 
are formulated for constructible ${\frak o}/\pi^r{\frak o}$-sheaves only, but using truncation with respect to the standard $t$-structure they extend to ${\frak o}/\pi^r{\frak o}$-adic complexes $K$ with bounded constructible
cohomology sheaves. Since $L_\chi$ is smooth on $X_R$ for any $N$-torsion character $\chi$, then  $$f: X_R \to S=Spec(R)$$ is also universally locally acyclic with respect to  $K_R^{(\pi)} \otimes_{\frak o} L_\chi$. Here $\chi$ is viewed
as a character $\chi: \pi_1(X,0) \to {\frak o}(\zeta_N)^* \subset Gl_{\frak o}({\frak o}[\zeta_N])$.
Now, if two of three complexes in a distinguished triangle are universally locally acyclic, then also the third is. This remark implies that all $K_V \otimes_{\frak o} L_\chi \otimes^{\rm L}_{\frak o} {\frak o}/\pi^r{\frak o}$ are universally locally acyclic for $r=1,2,...$ and our given ${\frak o}$-adic perverse sheaf $K\in {\cal B}$, if it is represented by the system $(K_r)_{r\geq 1}$ of perfect complexes $K_r$ such that  $K_r \cong K \otimes_{\frak o}^{\rm L} {\frak o}/\pi^r{\frak o}$.

\medskip
Now we may also consider $K$ as an object in the category $D_c^b(X,Q)$
or in the category $D_c^b(\overline X,Q)$ for $Q=Quot({\frak o})$. No matter in which way
$K\in {\cal B} \Longrightarrow \chi(K)=0$ [notice $\chi(K)^2=\chi(H^\bullet(\delta_B)) \cdot \chi(K) + \chi(T)=0$, since $\chi(T)=0$ and $\chi(H^\bullet(\delta_B))=0$.] Furthermore $K\in {\cal B}$ implies $K \cong {\cal P}_K$ and $\nu_K=\dim(A)$. [Indeed $K^\vee * K \cong K*K \cong H^\bullet(\delta_B) \cdot K \oplus T$ implies that ${\cal P}_K$ is either a summand of $T$ or as an indecomposable constituent of $H^\bullet(\delta_B)) \cdot K$ isomorphic to $K$.
In the first case $K$ is translation invariant under $X$, which is impossible by property 5).
Therefore $K\cong {\cal P}_K$. Since ${\cal P}_K[\pm\nu_K] \hookrightarrow H^\bullet(\delta_B)\cdot K$, this implies $\nu_K \leq \dim(A)$. On the other hand $\mu(K^\vee *K) = 0$ and $K[-\dim(A)] \hookrightarrow K^\vee *K$ imply 
$\nu_K=\mu(K) \geq \dim(A)$. Hence $\nu_K=\dim(A)$.]

\medskip 
For characters $\chi=(\chi_1,\chi_2)$ of $\pi_1(X,0)=\pi_1(A,0)\times \pi_1(B,0)$
notice that over the base field $\C$ we know that for almost all characters $\chi_1$
of $\pi_1(A,0)$ the direct image $Rp_*(K_{\chi_1}) \in D_c^b(B,Q)$ must be a locally free translation invariant perverse sheaf on $B$ of rank $r$. For fixed $\chi_1$ therefore
${\cal H}^\bullet(Rf_*(K_\chi))= H^\bullet(B,Rp_*(K_{\chi_1})_{\chi_2})$, considered 
with coefficients in ${\frak o}$, is a ${\frak o}$-torsion module for almost all characters 
$\chi_2$ of $\pi_1(B,0)$. Now for $K\in D_c^b(X,{\frak o})$ this information is entirely encoded
in the exact sequences for $r\to \infty$
$$0 \to {\cal H}^n(Rf_*(K_\chi))/\pi^r \to {\cal H}^n(Rf_*(K_\chi \otimes^{\rm L} {\frak o}/\pi^r{\frak o})) \to {\cal H}^{n+1}(Rf_*(K_\chi))[\pi^r] \to 0 \  $$
in the form that $ {\cal H}^\bullet(Rf_*(K_\chi \otimes^{\rm L} {\frak o}/\pi^r{\frak o}))$ has bounded
length independent from $r$.
Now consider the reduction $(\overline X,\overline K)$ of $(X,K)$, 
defined over the algebraic closure $\overline\kappa$ of the finite residue field $\kappa$ of $R$
with respect to the maximal ideal of $V$,
and the base change ring homomorphisms $V \hookrightarrow \C$ and $V \to \overline\kappa$. 
By the universal local acyclicity of the structure morphism $f$ for all $K_V \otimes_{\frak o} L_\chi \otimes^{\rm L}_{\frak o} {\frak o}/\pi^r{\frak o}$ the above bounded length conditions over $R$ or $V$  are inherited to the reduction, so  that  
$ {\cal H}^\bullet(R\overline f_*(\overline K_\chi \otimes^{\rm L} {\frak o}/\pi^r{\frak o}))$ has bounded
length independent from $r$, again  for
almost all characters $\chi_2$ of $\pi_1(\overline B,0) \cong \pi_1(B,0)$ with respect a fixed but arbitrary $\chi_1$ outside some finite set of exceptional characters $\chi_1$.   
The short exact
sequences 
$$0 \to {\cal H}^n(R\overline f_*(\overline  K_\chi))/\pi^r \to {\cal H}^n(R \overline  f_*(\overline  K_\chi \otimes^{\rm L} {\frak o}/\pi^r{\frak o})) \to {\cal H}^{n+1}(R\overline  f_*(\overline  K_\chi))[\pi^r] \to 0 \  $$
%
therefore imply that the $Q$-adic cohomology groups ${\cal H}^n(R\overline f_*(\overline K_\chi))$ vanish for almost all characters $\chi_2$ (with respect to the fixed $\chi_1$). Passing to
the algebraic closure $\Lambda$ of $Q$ we can apply the decomposition theorem
and obtain ${\cal H}^k(R\overline f_*(\overline K_\chi)) = \bigoplus_{i+j=k} H^i(\overline B, 
{}^p H^j (R\overline p_*(\overline K_\chi))$, and hence the $\Lambda$-adic perverse sheaves ${}^p H^j (R\overline p_*(\overline K_\chi))$ are acyclic for almost all $\chi_2$ (for fixed $\chi_1$). Since the perverse sheaves ${}^p H^j (R\overline p_*(\overline K_\chi))$ are pure $\Lambda$-adic Weil sheaves on $\overline B$,  by [W3, lemma 13] then all  ${}^p H^j (R\overline p_*(\overline K_\chi))$ are translation invariant perverse sheaves on $\overline B$. Hence $L=R\overline p_*(\overline K_\chi) =\bigoplus L_i[-i]$  for certain translation invariant perverse sheaves $L_i$. Now $L*L \cong H^\bullet(\delta_B) \cdot L \oplus R\overline p_*(\overline T_\chi)= H^\bullet(\delta_B) \cdot L $ for almost all $\chi_1$ implies $L=L_0$, since
$L_i[-i]*L_i[-i]$ contains nontrivial perverse sheaves in degree $-2i - \dim(B)= -2i-\dim(A)$
for any translation invariant nontrivial perverse sheaf $L_i$ on $\overline B$. 
Thus for fixed $\chi_1$ (for almost all $\chi_1$) the direct image $R\overline p_*(\overline K_\chi)$ itself is
a translation invariant perverse sheaf on $\overline B$. Up to a sign $(-1)^{\dim(B)}$ its rank  
is $\chi(L) = \chi(\overline K_\chi\vert_{\overline p^{-1}(0)}) =  
   \chi(\overline K\vert_{\overline p^{-1}(0)}) = \chi( R\overline p_*(\overline K)) = (-1)^{\dim(B)}\cdot r $. Using that $K$ and $L$ satisfy $K^\vee\cong K, L\cong L^\vee$ and
$L*L \cong H^\bullet(\delta_B)\cdot L$ and $\overline K_\chi*\overline K_\chi \cong H^\bullet(\delta_B) \oplus \overline T_\chi$,
which follows from the commutative diagram ($\bf *$) above, we are in the situation of proposition \ref{product} from the appendix. Hence this proposition implies
$$M^d(b)\cong \delta_0 \quad \quad \mbox{ for } b=0 \ .$$ In particular 
the generic rank $r$ of $L$ therefore is one. Hence $L$  is isomorphic to a direct sum of pairwise non-isomorphic translation invariant perverse sheaves of generic rank one for almost all torsion characters $\chi_1$ of $\pi_1(\overline A,0)$ (as explained at the beginning of this section). 
\qed

\medskip
From the last theorem we get

\begin{Proposition} \label{ToDo} $F(X)=\emptyset$ for complex abelian varieties $X$ isogenous to $A_1\times A_2$ with simple factors $A_1,A_2$ of dimension $\dim(A_1)=\dim(A_2)$.
\end{Proposition}
 
{\it Proof}. Assume there exists $K\in F(X)$. By theorem \ref{cp} we conclude that $M=K\vert p^{-1}(0)$ has Euler characteristic
one. This holds also for $K$ replaced by $K_\chi$ (for all $\chi$). 
Hence except for finitely many $\chi$  from $H^\bullet(A,M) \cong Rp_*(K_\chi)_0 $ and the fact that Euler characteristics of perverse sheaves are nonnegative [FK] we get 
that the perverse sheaf $M^d$ has Euler characteristic one
and $M^0,..,M^{d-1}$ are acyclic (using lemma \ref{Kill1}). Then all Jordan-H\"older constituents of
 $M^d$ are acyclic except for one, the perverse sheaf $\delta_0$ 
(lemma \ref{Kill1}) arising from the stalk spectral sequence analogous to lemma \ref{notzero}.
In fact this is also clear from an abstract point of view; the constituent
with Euler characteristic one is invertible in the Tannakian sense and therefore is a skyscraper sheaf
[KrW, prop.21b)]. This Tannakian argument carries over to
all fibers $F(b)=p^{-1}(b)$ and defines a unique perverse skyscraper  Jordan-H\"older constituent
in $M(b) = K_\chi\vert_{F(b)}$; in fact a perverse quotient sheaf of $M^d(b)$. Hence for every $b\in B$ this defines a point $x\in X$ with $p(x)=b$. The supports
of these skyscraper sheaves define a constructible set; its closure $S$ 
defines a birational morphism $p:S\to B$. On an open dense subset $U\subset S$
the morphism $p$ defines an isomorphism onto a dense open subet $V\subset B$.
By Milne [M, cor.3.6] the morphism $V \to U \to X$ extends to a homomorphism 
$s:B\to X$ (up to a translation). 

\medskip
For generic $b\in B$  we have an epimorphism of perverse sheaves
$$   M^d(b) \ \longrightarrow  \ \delta_{s(b)}  $$
so that $\delta_{s(b)}$ is the maximal perverse quotient of $M(b)$ with generic support
in a subvariety of dimension zero ([KW, lemma III.4.3]). 
The kernel is acyclic and it is nontrivial! 
[If it were trivial for generic $b$, then 
${\cal H}^{-d}(M^d) $ vanishes. This implies ${\cal H}^{-2d}(K\vert_{F(b)})=0$
for generic $b$. Hence ${\cal H}^{-\dim(X)}(K)$ vanishes at the generic point of $X$.
In fact the support is contained in $s(B)$. But this is impossible, since then
$K$ has support in a simple abelian subvariety and therefore is in $F(X)$.]
For generic $\chi$ the same conclusion also holds when $K$ is replaced by $K_\chi$ . 
Thus we always find acyclic nontrivial perverse subsheaves of $M^d$, which are
therefore $A$-invariant. In particular this implies $supp(K)=X$ (so in fact it would have not been necessary to suppose this). Now we apply verbatim the arguments of step 3) and step 5) of the proof of proposition \ref{fantom} and 
conclude as in proposition \ref{fantom} that the irreducible perverse sheaf $K$ must be translation invariant with respect to $A$; indeed $Kern(p)=A$ is simple and $supp(K)=X$. A contradiction.  \qed

\section{Appendix}

\bigskip\noindent
In this appendix we consider perverse sheaves on abelian varieties over finite fields and we use the Fourier transform as defined in [W3, section 7]
to discuss monoidal perverse sheaves $K_0$  on products of simple abelian varieties $X_0$ over a finite field $\kappa$. 
For a fixed embedding of $\kappa$ into an algebraic closure $k$, let $K$ resp. $X$ denote
the base extensions of $K_0$ resp. $X_0$ to the base $k$.  
In the following we regularly use notation from [W3, p.32 ff], often without further mentioning. 

\medskip
In this appendix we make the following assumptions:

\begin{itemize}
\item (A1) For simple abelian varieties $X$ over $k$ all perverse Weil sheaves $K$ on $X$, defined over a finite subfield $\kappa$ of $k$, that are monoidal in the sense of [W2], are translation invariant. 
\item (A2) If $K$ is an acyclic perverse Weil sheaf on a simple abelian variety $X$, for a finite subfield of $k$, then each irreducible perverse constituent  
of $K$ is acyclic (later to be applied for the perverse sheaves denoted $M^i$).
\end{itemize}

These assumptions are motivated by the fact that the analogous assumptions hold for simple abelian varieties over an algebraically closed field $k$ of characteristic zero,
as shown in [KrW]. We expect them to hold also in positive characteristic.

\medskip
Under the assumptions (A1) and (A2),  in corollary \ref{lc} we will 
show that any irreducible perverse sheaf $K$ defined over $\kappa$ on an abelian variety $X$ over $k$ with Euler characteristic zero is translation invariant with respect to some nontrivial abelian subvariety of $X$.  By the argument of proposition \ref{ToDo}, to show this it suffices to consider the case of irreducible monoidal perverse sheaves $K$ with Euler characteristic $\chi(K)=0$ on products $X$ of two simple abelian varieties $A,B$ over $k$. In addition, as shown in section \ref{next}, for the proof one can furthermore impose rather strict conditions on these monoidal perverse sheaves $K$ which naturally leads to a bunch of conditions that are formulated before the next proposition \ref{product}. Then, using our arguments from section \ref{specialsheaves}, the desired corollary \ref{lc} will be a consequence of proposition \ref{product}. The assumptions (A1) and (A2) are needed to ensure that the technical conditions imposed in proposition \ref{product} for the monoidal perverse sheaf $K$ do hold. 

\medskip
{\it Remark}. In the proof of theorem \ref{cp} we also applied proposition \ref{product}. However there, by the reduction procedure, the required assumptions for proposition \ref{product} are inherited from the characteristic zero situation where (A1) and (A2) hold unconditionally.  Hence for the proof of theorem \ref{cp},   fortunately, it is  not necessary to postulate the conditions (A1) and (A2)   
from above.

\medskip
{\it Notations}. Assume $\Lambda= \overline\Q_l$ for some prime $l$ different from $p$.
Let $\kappa$ be a finite field of characteristic $p$ with $q$ elements and $\kappa_m$ be a finite extension field of $\kappa$ of degree $m$ with $q^m$ elements. Let $F=F_{\kappa}$ denote the Frobenius automorphism of $\kappa$ so that $F_{\kappa_m}=F^m$. Let $X_0$ be an abelian variety over $\kappa$ and $X$ its scalar extension  to the algebraic closure $k$ of $\kappa$. Let  $T_0$ in $D_c^b(X_0,\Lambda)$ denote a semisimple translation invariant complex with the property ${}^p H^\nu(T_0)=0$ for $\nu < \dim(X)$ and $\nu >-\dim(X)$. Let $H^\bullet_0 = \oplus_{i=-d}^{d} H_0^i[-i]$ denote some fixed complex in $D_c^b(Spec(\kappa),\Lambda)$ and $ H^\bullet$ the associated $F$-module and $h_m:= Tr(F^m;H^\bullet)$. We assume $\chi(H^\bullet)=\sum_\nu (-1)^\nu \dim_\Lambda(H^\nu)=0$. Concerning this, we remark that later we apply this for $X_0 = A_0 \times B_0$ and $H^\bullet = H^\bullet(B,\delta_B)$ considered as a $F_\kappa$-module, a situation related to  section \ref{next}. Let $K_0\in Perv(X_0,\Lambda)$ be a perverse sheaf, whose pullback $K$ to $X$ is a simple perverse sheaf on $X$.  
Let ${\cal S}(K)$ denote the spectrum of $K$, i.e. the set of \lq{continuous}\rq\ characters $\chi: \pi_1(X,0) \to \Lambda^*$ of the etale fundamental group of $X$ for which $H^\bullet(X,K_\chi)\neq H^0(X,K_\chi)$ holds; in the following it suffices  to consider torsion characters $\chi$. For $T$ we similarly define the spectrum ${\cal S}(T)$ of $T$ by considering the perverse cohomology sheaves ${}^p H^i(T)$ of $T$. Let $X_m$ denote the $\kappa_m$-rational points $X(\kappa_m)$ of $X$ using some fixed embeddings $\kappa_m\hookrightarrow k$. Similarly $A_m=A(\kappa_m)$ and $B_m=B(\kappa_m)$. Let $X_m^*$ be the group of characters $\chi: X_m \to \Lambda^*$, also considered as a group of characters of $\pi_1(X,0)$ via Lang torsors as in [W3, section 7]. This  defines the points of the spectrum
${\cal S}_m(K) ={\cal S}(K) \cap X_m^*$ and ditto ${\cal S}_m(T)$ as the union $\bigcup_i {\cal S}({}^p H^i(T))\cap X_m^*$. 
For a (super)natural number $N$ let $X[N^\infty]$ denote the set of points in $X(k)$ annihilated by some power of $N$.

\goodbreak

\medskip
Define $\Lambda_{0}^* = \{ \alpha\in \Lambda^* \ \vert \ \vert \alpha\vert_{\rm v} = 1 $ for all archimedean and nonarchimedean valuation $\vert . \vert_{\rm v}$ not over $p \}$.  If there exists integers $n$ depending on $\alpha$ with $\alpha q^{n/2} \in \Lambda_{0}^*$, then by definition  $\alpha\in \Lambda_{mot}^*$ and we write $w(\alpha)=-n$. 
Following [Dr2], a weakly motivic Weil complex $K$ on $X$ is a complex with the property that for all cohomology sheaves ${\cal H}^\nu(K)$ the eigenvalues of the Frobenius $Fr_{x}$ on the stalks ${\cal H}^\nu(K)_{\overline x}$ at geometric points $\overline x$ over closed points $x$ of $X$ have algebraic eigenvalues $\alpha\in \Lambda_{mot}^*$, i.e. there exists integers $n$ depending on $\alpha$ with $\alpha q_x^{n/2} \in \Lambda_{0}^*$. By [Dr2, thm.B.3] the weakly motivic complexes define a thick triangulated full subcategory of $D_c^b(X_0,\Lambda)$. For $k$-morphisms $f:X\to Y$ these are stable under the functors
$f_!,f_*,f^*,f^!,\boxtimes,D,R{\cal H}om$. For any weakly motivic complex $K$ the functions
$f^K_m(x) = \sum_\nu (-1)^\nu Tr(F_\kappa^m; {\cal H}^\nu(K)_{\overline x})$ have values in a number field (independent from $m$) and these values are integral over $\Z[p^{-1}]$. 

\medskip
As in the next proposition, now assume for the irreducible perverse sheaf $K$ 
$$  K^\vee \cong K \quad \mbox{ and } \quad K^\vee*K \cong H^\bullet\cdot K \oplus T \ .$$
As shown in proposition \ref{deco}, for monoidal perverse sheaves $K$ on a product $X$ of simple abelian varieties  this automatically holds as a consequence of our assumption (A1) on simple abelian varieties.
By replacing  the irreducible Weil complex $K$ by a generalized Tate twist $K(\alpha)$ for some $\alpha\in \Lambda^*$, we may achieve that the determinant $det(E)$ of the smooth local coefficient system $E$, defining $K$ on some open dense subset of its support, has finite order [L]. Then $K_0$ is weakly motivic (see [Dr2]). As an immediate consequence also
$H_0^\bullet$ and $T_0$, defined by 
$$ K_0*K_0 \cong H^\bullet_0\cdot K_0 \oplus T_0 \ ,$$ 
are weakly motivic. Notice $K_0^\vee \cong K_0(\alpha)$. By the determinant condition $\alpha$ is a root of unity. Hence $K$ is pure of weight $0$. By a suitable finite extension of the finite base field $\kappa$ we may then assume $K_0^\vee \cong K_0$ so that
$$  K_0^\vee \cong K_0 \quad \mbox{ and } \quad K_0*K_0 \cong H^\bullet_0\cdot K_0 \oplus T_0 \ .$$ 
Instead for $(K_0,H_0^\bullet,T_0)$ this also holds for the triples
$(K_{0\chi},H_0^\bullet,T_{0\chi})$  obtained from twisting with a  character $\chi\in X_m^*$. Since the Euler characteristic
of $H^\bullet$ and $T$ is zero, the Euler characteristic of all $K_\chi$ is zero. This implies $\chi\in {\cal S}(K)$ iff $H^\bullet(X,K_\chi)\neq 0$. In particular $\widehat f^K_m(\chi)=0$ holds for all $\chi\notin {\cal S}(K)$ by the Grothendieck-Lefschetz trace formula (see [W3], section 6).

\medskip
If $\chi\notin {\cal S}(T)$, then $H^\bullet(X,K_\chi)^{\otimes 2} \cong H^\bullet \otimes_\Lambda H^\bullet(X,K_\chi)$. Hence $\widehat f^K_m(\chi)\neq 0$
implies 
$$  \widehat f^K_m(\chi)= Tr(F^m; H^\bullet(X,K_\chi)) =:  h_m \ $$
for $h_m = Tr(F^m,H^\bullet)$.
We remark that $U\otimes_\Lambda W \cong V \otimes_\Lambda W$ and $W\neq 0$ for semisimple $F$-modules $U,V,W$ implies $U\cong V$ after a finite field extension. Hence for the characteristic function $1_{{\cal S}_m(K)}$ of ${\cal S}_m(K)$ more precisely $$ \widehat f^K_m(\chi) \ = \ h_m \cdot 1_{{\cal S}_m(K)}(\chi) \ + \ \mbox{function supported in }\ {\cal S}(T) $$ holds for all $m$, after replacing $\kappa$ by a suitable finite field extension. Hint:  If the $F$-eigenvalue multi sets $A,B,C$ of $U,V,W$ satisfy $A + C = B + C$ in $C^*\otimes_{\Z}\R$, to show $A=B$ (as sets with multiplicities) simply use induction on $\dim(U)=\dim(V)$. For the
induction step choose some arbitrary (lexicographic) ordering on the finite dimensional $\R$-vector space generated in $C^*\otimes_{\Z}\R$ by the elements in $A,B,C$ that respects addition. Then compare $a+c$ and $b+c$ for the largest elements $a,b,c$ in $A,B,C$ to conclude $a=b$.

\medskip
Now let
$X_0=A_0\times B_0$ be a product of two simple abelian varieties $A_0$ and $B_0$ defined over $\kappa$. Let $p_0: X_0 \to B_0$ be the projection onto the second factor.
The characters $\chi \in X_m^*$ correspond to pairs $(\chi_1,\chi_2)$ of characters $\chi_1\in A_m^*$ and $\chi_2\in B_m^* $; we assume that except for $\chi_1$ in a {\it finite set} $\Sigma$, for every torsion character $\chi_1$ there exist finitely many characters $\psi_i$ depending on $\chi_1$
such that $$ Rp_*(K_{\chi_1}) \cong \bigoplus_{i=1}^r \ m_i \cdot \delta^{\psi_i}_B \ .$$ 
Since $F(B)=\emptyset$ holds by assumption (A1), for monoidal perverse sheaves this follows from lemma \ref{strange}.  

\medskip
In this situation, by section \ref{specialsheaves} the multiplicities $m_i$ are zero or one, and the number $r=r(\chi_1)$ of characters that appear
with multiplicity $m_i >0$ is independent of $\chi_1$ and  is called $rank(Rp_*(K))$.  Up to a sign this rank  is the Euler characteristic
$r = (-1)^d \chi(M(b))$ of the complex $M(b)=K\vert_{A\times \{b\}}$ on ${A\times \{b\}}$ (for $b\in B$). Attached to the sheaf complex $M:=M(0)$ on $A\times \{0\}$ consider
its perverse sheaves $M^i ={}^p H^{-i}(M(0))$.  Except for $i=0,..,d=\dim(B)$ the $M^i$ are zero. For $\chi_1\notin \Sigma$ the $M^i_\chi$ are acyclic on $A = A\times \{ 0\}$ for $i\neq d$, and similarly for 
$M^i(b)$. By our second assumption (A2) all simple perverse constituents of these $M^i, i\neq d$ are
acyclic and hence by the first assumption (A1) the $M^i,i\neq d$ then are translation invariant perverse sheaves on $A$. Since $K$ is pure of weight $0$, this implies $w(M)\leq 0$ and hence
$w(M^i) \leq -i $ for all $i$. 

\goodbreak

\begin{Proposition} \label{product}
For $X_0=A_0\times B_0$ and simple abelian varieties $A_0,B_0$ of dimension $d=\dim(A_0)=\dim(B_0)$ let $K_0$ be a simple monoidal perverse sheaf on $X_0$ with $\nu_K=d$ such that $K^\vee \cong K$ and
$K^\vee*K \cong H^\bullet \cdot K \oplus T$ holds for a translation invariant complex $T=\bigoplus_{i=a}^b T_i[-i]$ on $X$ with  $-\dim(X)< a  \leq b<\dim(X)$  and a graded $\Lambda$-vectorspace\footnote{$\chi(H^\bullet)=0$ holds for simple monoidal perverse sheaves $K_0$ (see [W2]).}  and graded $F$-module $H^\bullet = \bigoplus_{i=-d}^d H^i[-i]$ such that $H^{\pm d} \cong \Lambda$. Suppose that the support of $K$ is not contained in a fiber of the projection $p_0: X_0\to B_0$ and that ${\cal H}^{-d}(K)=\delta_0$.  Assume there exists a finite set $\Sigma$ such that for all characters $\chi_1 \notin \Sigma$ of $\pi_1(A,0)$ the direct image $Rp_*(K_{\chi_1})$ is a perverse translation invariant sheaf on $B$. Assume that all simple perverse constituents of the associated perverse sheaves $M^i, i\neq d$ are translation invariant perverse sheaves on $A$. Then $$M^d \cong \delta_0$$ holds as perverse sheaves on $p^{-1}(0)=A\times \{ 0\}$ and in particular the rank $r$ of the locally constant sheaf $Rp_*(K)$ is $$  rank(Rp_*(K)) = 1 \ .$$
\end{Proposition}

\medskip
{\it Proof}. Suppose $\chi=(\chi_1,\chi_2)\notin {\cal S}(T)$ and $\chi_1\notin \Sigma$. This excludes finitely many characters $\chi_1$ of $\pi_1(A,0)$. By our assumptions then  $$Rp_*(K_\chi)= \bigoplus_{i=1}^r \ \delta_B^{\psi_i}$$ holds with a multiplicity free
collection of $r$ characters $\psi_i$ depending on $\chi$ (a special case of [W2,prop. 2]). It suffices to consider $\chi=(\chi_1,1)$ for torsion characters $\chi_1$,  outside a finite set of exceptional characters.

\medskip Step 1)
As explained in [W3, section 7], under the above assumptions on $\chi\in X_m^*$ we obtain an isomorphism of perverse Weil sheaves over $\kappa_m$ $$Rp_{0*}(K_{0\chi}) \cong \bigoplus_j \ Ind_{\kappa_{mr_j}}^{\kappa_m}(\delta^{\psi_j}_{B_0})(\alpha_j) $$  for some $\alpha_j, r_j$ and $\psi_j$ depending on $\chi$. Since $K_{0\chi}$ is pure of weight 0, 
$Rp_{0*}(K_{0\chi})$ is pure of weight zero. Hence $w(\alpha_j)=0$ holds for the weight for all $j$.

\medskip
Step 2) For $\chi\in A_m^*$, $b\in B_m$ let $\{\psi_1,..,\psi_s\}$ be the subset  of characters
$\psi_j$ from step 1) with the property $r_j=1$ (i.e. $\psi_j \in B_m^*$). Then  
example 2 and 3 of [W3, section 7] show  
$$ f^{Rp_*K_\chi}_m(b) \ = \ (-1)^d q_m^{-d/2} \cdot \sum_{j=1}^s  \ \alpha_j^m \cdot \psi_j(b)^{-1} \ .$$  
The number of summands $s=s_m(\chi)$ in this formula may depend on $\chi \in {\cal S}_m(K)$ and $m$. Also the twist factors $\alpha_j^m$ may a priori depend on $\chi \in {\cal S}_m(K)$.

\medskip
Step 3) Notice that a  twist of $Rp_*K_{\chi} $ with the inverse of $\psi_j = (1,\psi_j)$ for $j=1,..,s$ gives  $Rp_*K_{\chi/ \psi_j} = (Rp_*K_{\chi})_{\psi_j^{-1}} = \delta_B(\alpha_j) \oplus$ acyclic perverse
 sheaves on $B$.
The twist factors $\alpha_j\in \Lambda_{mot}^*$ for $1\leq j \leq s$ defining the \lq{Tate twists}\rq\ in the last formula can therefore (for all $m$) be computed by the nonvanishing number 
$$ \alpha_j^m \cdot Tr(F_{\kappa}^m; H^\bullet(B,\delta_B))   = Tr(F_{\kappa}^m; H^\bullet(B,\delta_B(\alpha_j))  $$
$$ = Tr(F_{\kappa}^m; H^\bullet(B, (Rp_*K_{\chi/\psi_j}))  = Tr(F_{\kappa}^m; H^\bullet(X, K_{\chi/\psi_j}))  = \widehat f^K_m(\chi/\psi_j) \ . $$
Now $\chi/\psi_j \in {\cal S}_m(K)$ follows from $\chi\in {\cal S}_m(K)$; furthermore $\chi/\psi_j \notin {\cal S}_m(T)$ holds, since we discarded a finite set of expectional characters $\chi=(\chi_1,1)$.
Hence using $\chi/\psi_j \in {\cal S}_m(K)$ we get  
$$ \widehat f^K_m(\chi/\psi_j)  = h_m = Tr(F_{\kappa}^m; H^\bullet) \ ,$$
as already shown.
Hence the $\alpha_j^m$ are independent from $j$ and independent from $\chi=\chi_1\notin \Sigma$ in ${\cal S}_m(K)$.
The arguments used in the proof of proposition \ref{deco}  also show that $Tr(F_{\kappa}^m; H^\bullet)$
is $\alpha^m \cdot Tr(F_{\kappa}^m; H^\bullet(B,\delta_B))$ for some
twist factor $\alpha$.  Indeed, by the proof of corollary \ref{funny},  the factor $H^\bullet$ is defined by the formulas $K*K \cong H^\bullet \otimes_\Lambda P \oplus T$ and
$Rc_*c^*[d](P) = H^\bullet \otimes_\Lambda P$, which are identities of Weil sheaves 
and $c$ over $k$ comes from a homomorphism $c_0: A_0\times B_0^2 \to A_0\times B_0$ defined over $\kappa$  with kernel isomorphic to $B_0$. Hence $H^\bullet = H^\bullet(B_0,\Lambda[d])$.  Furthermore $P\cong K$ over $k$ and hence $P_0\cong K_0(\beta)$  for some $\beta$.
This implies $H^\bullet = H^\bullet(B,\delta_B)(\alpha)$ for $\alpha=q^{d/2}\beta\in \Lambda_{mot}^*$.
Notice $Tr(F_{\kappa}^m; H^\bullet(B,\delta_B)) = (-1)^d q_m^{-d/2}\# B_m \neq 0$ for all
$m$, hence $\alpha_j^m=\alpha^m$.
In particular $w(\alpha)=0$ and $w(\beta)=-d$ holds for the weights, and our formula obtained in step 2) simplifies to 
$$ f^{Rp_*K_\chi}_m(b) \ = \ (-1)^d \beta^m \cdot \sum_{j=1}^s  \ \psi_j(b)^{-1} \ $$ 
for some fixed $\beta$ of weight $w(\beta)=-d$.

\medskip
Step 4) Since $K_\chi$ and $Rp_*(K_\chi)$ are weakly motivic, there exists a subring ${\frak o} \subset \Lambda$ finite over $\Z_l$ with prime ideal generated by $\pi$ such that $f_m^{K_\chi}(x)$ and also $f^{Rp_*(K_\chi)}_m(b)$ have values
in ${\frak o}$.  Notice ${\frak o}$ may depend on $\chi$.
For $\chi',\chi \in A_m^*$ and $\chi',\chi \notin {\cal S}(T)$ such that $\chi' = \chi\chi_l$
holds for an $l$-power torsion character $\chi_l$, then from step 2) and 3)  for all $b\in B_m$
we get the following equality in  the residue field ${\mathbb F} = {\frak o}/\pi$
$$ f_m^{Rp_*K_{\chi'}}(b) \ \equiv \   f_m^{Rp_*K_{\chi}}(b) \ .$$ 
This follows from $f_m^{Rp_*(K_{\chi'})}(b) = \sum_{a\in A_m} f_m^{K_{\chi'}}(a,b)
= \sum_{a\in A_m} f_m^{K}(a,b)\cdot \chi'(a,b)^{-1}$ and the property $f_m^K(a,b)\in {\frak o}$, 
using $\chi_l(a,b)\equiv 1$ mod $\pi$ and  
$\chi'(a,b)^{-1} \equiv \chi(a,b)^{-1}$ in $\mathbb F$. 

\medskip
Step 5) Writing $f_m^{Rp_*K_{\chi}}(b)$ as a finite linear combination
of $s$ characters $\psi_j: B_m \to {\mathbb F}^*$ as in step 3), 
for $b=0$ we obtain
$$  f_m^{Rp_*K_{\chi}}(0) \ = \ (-1)^d \cdot \beta^m \cdot s $$
where the number  of characters $s=s_m(\chi)$ may depend on $\chi$ and $m$.
By the congruence formula from step 4) 
for $\chi'=\chi\chi_l$  then
$$ \beta^m \cdot \sum_{i=1}^{s(\chi')} \psi'_j(b)^{-1} \equiv   \beta^m \cdot \sum_{i=1}^{s(\chi)} \psi_j(b)^{-1} \ $$
follows for all $b\in B_m$. But $\beta\in \Lambda_{mot}^*$ (so that $\beta \not\equiv 0$ in $\mathbb F$) and the linear independency of characters  $B_m \to \mathbb F^*$ then implies the congruence $s_m(\chi')\equiv s_m(\chi)$.  This congruence forces an equality $s_m(\chi')=s_m(\chi)$ if $l $ is larger than the Euler characteristic $r$
of $M(0)$. Indeed, $r$ is an upper bound for $s=s_m(\chi)$  independent from $\chi$ and $m$. Hence for $l>r$ we 
get an equality
$$  f_m^{Rp_*K_{\chi'}}(0) = (-1)^d \cdot \beta^m \cdot s_m(\chi')  = (-1)^d \cdot \beta^m \cdot s_m(\chi)  =  f_m^{Rp_*K_{\chi}}(0) \ $$
whenever $\chi' =\chi\chi_l$ satisfies the conditions of step 4).

\medskip
Step 6) By [W3, thm. 5] this can be applied for any primes $\ell' \neq p$ and $\ell' > r$ instead of the given prime $\ell$. 
Let $N= p \prod_{\ell'\leq r} \ell'$  be the product of the finitely many remaining primes.
Then  $\chi\mapsto f_m^{Rp_*(K_\chi)}(0)$ is for all $m$ an 
$A_m^*[l'^\infty]$-invariant function on $A_m^*$ for the primes $\ell' \not\vert N$, at least on the complement of the set ${\cal S}_m(T)$. For simplicity we can assume ${\cal S}(T) \subseteq {\cal S}_m(T)$ by enlarging $\kappa$.
Let $N'$ be the supernatural number $N'= \prod_{\ell' \not\in N} \ell'$. 
In fact, for any $m$ we can extend   this  function to  a constant function on each $A_m^*[N'^\infty]$-coset in $A_m^*$. This defines $A_m^*[N']$-invariant functions $\widehat f_m(\chi)$ on $A_m$ for each $m$.  
  
\medskip
Step 7) From the last step  we get for the extensions $\widehat f_m(\chi)$ defined in step 6) the following formula  
$$ \widehat f_m^M(\chi_1) = f_m^{Rp_*(K_\chi)}(0) \ =\ \widehat f^K_m(\chi) =   \widehat f_m(\chi) + \mbox{ function with support in } {\cal S}_m(T) \ $$
for $\chi=(\chi_1,1)$ and all $\chi_1\in A_m^*$ and all $m$.

\medskip
Step 8) The perverse cohomology sheaves $M^i$ of $M$ for $i\neq d$ are translation invariant on $A\times \{0 \}$.
Hence 
in the formula of step 7) we can replace $M$ by $M^d$  after replacing ${\cal S}(T)$ by a suitable larger finite set ${\cal S}$ of characters (that is independent of $m$).  
By Fourier inversion then for all $a\in A_m$ the formula
$$ f_m^{M^d}(a) \ =\ f_m(a) \ + \ g_m(a) \quad , \quad
 g_m(a) \ =\ \sum_{\chi\in {\cal S}(T)} \ g(m,\chi)\cdot f_m^{\delta_A^\chi(\alpha_\chi)}(a) $$
holds for the Fourier inverse $f_m(a)$ of $\widehat f_m(\chi)$. Notice that by step 6) for all $m$ the functions $f_m$ have support
$$ supp(f_m) \subseteq A_m[N^\infty]\ .$$ 
By enlarging $\kappa$ we may assume  that all characters in ${\cal S}$
are in $A(\kappa)^*$, and by enlarging $N$ we may assume also that $\chi^N=1$ holds for $\chi\in {\cal S}$. 
This implies $g_m(a + x)=g_m(a)$ for all $x$ in some subgroup $U_m$ of $A(\kappa_m)$ with
index $[A(\kappa_m):U_m]$ dividing $N^{2\dim(X)}$. In particular, this holds for all points $x\in A(\kappa_m)[l^\infty]$ for a prime $l$ not dividing $N$. Let us fix such a prime $l$.

\medskip
Step 9) The semisimplification $(M^d)^{ss}$ of the perverse sheaf $M^d$ on the simple abelian variety $A$ can be written as a sum
$(M^d)^{ss} = T \oplus \delta$, where $T$ is translation 
invariant and where $\delta$ is clean in the sense that it does not contain translation invariant summands. Let $Y$ be an irreducible component of the support $supp(\delta)$ of $\delta$ whose dimension is maximal among in $supp(\delta)$. If $\dim(Y) < \dim(A)$, choose a Zariski open  subset $Y'\subset Y$ which is disjoint to other irreducible components of $supp(\delta)$ and
so that $\delta\vert_{Y'}$ is smooth, i.e. a smooth sheaf up to a complex shift. The Zariski closure of  $A(k)[l^\infty]$ in $A$ is $A$.
For $\delta\neq 0$ there exists an
integer $m$ and a point $x\in A(\kappa_m)[l^\infty]$ such that $Y'\cap (-x +supp(\delta))$
is a proper subset of $Y'$. If we replace $Y'$ by some open dense subset, we may therefore
assume $(Y' +x)\cap supp(\delta)=\emptyset$. 
Therefore we can find $m_0,a\in Y'(\kappa_{m_0})$ so that 
 $f_m^\delta(a+x)=0$ holds for all $m$ and so that $f_m^\delta(a)\neq 0$ holds for infinitely many $m$.
In the Grothendieck group of perverse sheaves on $A$ then the linear combination 
$$  T_x^*\bigl((M^d)^{ss}\bigr)\ - \ (M^d)^{ss} \ = \ T_x^*(\delta)\ - \ \delta $$
is nontrivial, with certain components having generic support of dimension $\dim(Y)$ not cancelling away. On the other hand by step 8), for varying $m$ the associated functions of this element in the Grothendieck group are $a\mapsto f_m^{M^d}(a + x)
- f^{M^d}(a) = f_m(a + x) + g_m(a+x) - f_m(a) - g_m(a) = f_m(a+x) - f_m(a)$, and their support
hence is contained in $A(\kappa_m)[(lN)^\infty]$. In view of the first assertion of the remark following [W3, cor. 3], this contradicts  [W3, lemma 19] 
(applied for $lN$ instead of $N$) unless $\dim(Y)\leq 0$. So, for $\dim(Y)\neq \dim(A)$
the semisimplification $(M^d)^{ss}$ of the perverse sheaf $M^d$ on $A$ is $$ (M^d)^{ss} \ \cong \ T \oplus \delta \ $$
for some perverse Weil sheaf $\delta$ on $A$ with support of dimension 0 
and some translation invariant perverse Weil sheaf $T$ on $A$.
 
\medskip
Step 10) For $\delta\neq 0$, we claim that $\dim(Y)=\dim(A)$ is impossible. Otherwise
on some open dense subset $Y'\subset A$ the semisimple perverse Weil sheaf $\delta$ becomes smooth $\delta\vert_{Y'} = E[d]$ for some smooth $l$-adic sheaf $E\neq 0$. Since none of the irreducible smooth sheaf components $E_i$ of $E\cong \bigoplus_i m_i \cdot E_i$ is translation invariant on $X$, we can choose $x$ as in step 9) so that
$\delta$ and $T_x^*(\delta)$ have no common simple smooth sheaf constituent on $Y'$. Indeed, the stabilizer in $X$ of each constituent is a finite group.  
Furthermore, as in the last step  
for all $m$ the functions
$$ f^{T^*_x(\delta)}_m(a) - f^{\delta}_m(a) $$
have support in $A_m[(lN)^\infty]$. These two mentioned properties of $\delta$ are obviously stable under (generalized) Tate twists; so for simplicity we may now temporarily assume that the maximal weight of a nontrivial
simple constituent of $\delta$ is $w=0$ resp. of $E$ is $-d$.  
For the restrictions of $\delta$ etc. to $Y'$ (and by abuse of notation we also write $\delta$ etc.) consider the corresponding functions $f_m^\delta$ etc. on $Y'\cap X_m$.
For $m\to \infty$ the left side of the scalar products on $Y\cap X_m$
 $$ (f_m^\delta, f_m^\delta - f_m^{T_x^*(\delta)}) = (f_m^\delta, f_m^\delta) - (f_m^\delta, f_m^{T_x^*(\delta)}) \ $$
can then be estimated from above by $2\# (supp(f_m^\delta - f_m^{T_x^*(\delta)}) \cdot (q_m^{w(E)/2})^2 \cdot 
rank(E)^2$.  Hence 
$2\#{A_m[(lN)^\infty]}\cdot q_m^{-d} \cdot rank(E)^2$ is an upper bound. We may enlarge $\kappa$
so that $A_0(\kappa)$ contains a torsion point of order $M$ for some integer $M$ prime to $N\cdot l$
such that $2\cdot rank(E)^2 < M \cdot \sum_i m_i^2$.  Then the left side becomes $<  \sum_i m_i^2$
for large $m$, since  $\frac{\#{A_m[(lN)^\infty]}}{\#A_m}\leq \frac{1}{M}$ and  $\#{A_m}q_m^{-d}$ converges to one. By Cebotarev, $(f_m^\delta, f_m^{T_x^*(\delta)})$ on the right side converges to zero ([W3, cor. 2]), whereas for any $\varepsilon >0$ there exist infinitely many integers $m$ such that $(f_m^\delta, f_m^\delta) > \sum_i m_i^2 - \varepsilon$ holds, using the second assertion of the remark following [W3, cor. 3]. This gives a contradiction, so $\dim(Y)< \dim(A)$.  
Hence, the support of $\delta$ on $A$ has dimension 0 or $\delta=0$.  

\medskip
Step 11)  By the last step $\delta \cong \bigoplus_i m_i \cdot {\delta_{x_i}}(\beta_i)$, for finitely many closed points $x_i\in A$ with certain multiplicities $m_i\geq 1$ and twists $\beta_i \in \Lambda_{mot}^*$ (unless $\delta=0$). We choose $\kappa$ large enough so that the points $x_i$ are
$\kappa$-rational points.

\medskip
Step 12) From step 9), 10), 11) for all $m$ and for almost all characters $\chi=\chi_1$ we obtain
$$ f_m^{Rp_*(K_\chi)}(0)= \sum_{a\in A_m} f_m^{M_\chi}(a) = \sum_{a\in A_m}
f_m^{M^d_\chi}(a) = \sum_{a\in A_m} f_m^{\delta_\chi}(a) = \sum_i m_i \cdot \beta_i^m \cdot \chi(x_i)^{-1} 
\ .$$ 

\medskip
Step 13) For pointwise weights we have
$w(M^d) \leq w(M) - d \leq w(K) - d = -d$. We claim that the weights of the stalks of $\delta$ at the points $x_i$ in the support are equal to $-d$;  as a skyscraper sheaf then $\delta$ is pure of weight $-d$. Indeed,  
by comparing step 11) and step 5) for all $m$ and all $\chi \in A_m^*$   
(with finitely many exceptions),  we get 
$$ (-1)^d \cdot s_m(\chi) = \sum_i \ m_i \cdot \theta_i^m \cdot \chi(x_i)^{-1}  \  $$
for the numbers $\theta_i= \beta_i/\beta$ in $\Lambda_{mot}^*$. Hence,
by regrouping into pairwise different $\theta_k$, this shows 
$$ (-1)^d \cdot s_m(\chi) = \sum_k \ \theta_k^m \cdot m_k(\chi) \quad , \quad  m_k(\chi) = \sum_i \ m_{ki} \cdot \chi(x_i)^{-1} \neq 0 \  $$
with integer coefficients $m_{ki}\geq 0$.  
Since $w(\beta) = -d$ and $w(\beta_i)\leq -d$ and $s_m(\chi)$ is an integer, if $w(\beta_i)< -d$ would hold for one of the $\beta_i$, this would give a contradiction for $m\to \infty$; in fact, it would imply $w(\theta_k)<0$ 
for some k. Then choose $\chi\in A_m^*$  such that $\chi(x_i)=1$ holds for all $i$ 
and for some $m=m_0$. We can also assume $s_{m_0}(\chi)=r$, by enlarging $m_0$ if necessary.  Then, for $m=\nu m_0$ and $\nu\to \infty$ we have $s_m(\chi)=r$, and the formulas above 
give a contradiction. This proves $w(\beta_i)=-d$ for all $i$. 

\medskip
Step 14)  Since $w(M^d)\leq -d$ and since $\delta$ is pure of highest weight $-d$ as shown in step 13), 
the weight filtration on $M^d$ and the decomposition theorem imply the existence 
of an exact sequence
$$ 0 \to T \to M^d \to \delta \to 0 $$
of perverse sheaves on $A$.

\medskip
Step 15)  This in turn implies $T=0$ and hence $M^d \cong \delta$. Indeed, the existence of a translation invariant perverse subsheaf $T\neq 0$ in $M^d$ would imply that $K$ is translation invariant under $A$, by the argument used in step 2 and 3 of the proof of proposition \ref{fantom} (note that the argument there does not depend on the choice of the base point $b$).  Hence by our assumption ${\cal H}^{-d}(K) \cong \delta_0$, the restriction $M$ of the perverse sheaf 
$K$ to $A\times \{ 0\}$ can not be translation invariant under $A$. Therefore 
$$M^d \cong \delta \neq 0\ .$$ 
\indent 
Step 16) $K$ is a simple monoidal perverse sheaf with ${\cal H}^{-d}(K)_0 \cong \Lambda$ by  assumption. The stalk spectral sequence of section 1 (before lemma \ref{notzero}) therefore
gives an exact sequence of etale sheaves on $A\times \{ 0\}$
$$ \xymatrix{ 0 \ar[r] & {\cal H}^{-d}(M^1) \ar[r]^{\partial_d} & {\cal H}^{0}(M^d) \ar[r] &
{\cal H}^{-d}(K)\vert_{A\times\{0\}} \ar[r] & {\cal H}^{-d}(M_0)  } \ .$$
Indeed, its higher differentials $\partial_i: {\cal H}^{-i}(M^j) \to {\cal H}^0(M^d)$ for 
$i+j=d+1$, $i=2,..,d-1$ vanish, because the perverse sheaves $M^i$ for $i < d$ are translation invariant on $A$ and so their cohomology sheaves are zero
except in degree $-d$.  Our assumptions that the support of $K_0$ is not contained in a fiber of $p_0$ implies $M^0=0$. Hence, the right term ${\cal H}^{-d}(M_0) $ vanishes.
Finally notice, $w(M^1)\leq -1$ and hence $w({\cal H}^{-d}(M^1)) \leq -d-1$, but $M^d$ is a pure skyscraper sheaf
of weight $w=-d$. So, by weight reasons we conclude $\partial_d=0$ and hence ${\cal H}^0(M^d)_a  \cong   {\cal H}^{-d}(K)_a $ holds for all $a\in A \times \{ 0\}$. Since $ {\cal H}^{-d}(K) = \delta_0 $ by our assumptions, therefore   
$$\delta  \ \cong \  {\cal H}^{-d}(K)\vert_{A\times\{0\}} = \delta_0 \ .$$
Since $\nu_K=d$, the support of ${\cal H}^{-d}(K)$ can be identified with the finite stabilizer group of the monoidal perverse sheaf $K$ and each stalk has $\Lambda$-dimension one by [W2, lemma 1, part 5 and 7].

\medskip
Step 17)
Since $\delta=\delta_0$, the Euler characteristic of $Rp_*(K)$ or $M$  is $(-1)^d$ using $(-1)^d \chi(M) = \chi(M^d)=
\chi(\delta_0)= 1$. Hence $rank(Rp_*(K)) = (-1)^d \cdot \chi(M) = 1$.
\qed

\bigskip
{\it Remark}. Applying proposition \ref{product} as in the proof of theorem \ref{cp}, the roles of $A$ and $B$ can be interchanged. Using this, then step 3 of our proof implies $\#A_m=c^m \cdot h_m=\# B_m$ for almost all $m$. Since $\sum_{m\geq m_0}^\infty \# A_m t^m$ determines the $L$-function of
$A$, and similar for $B$, we get $Tr(F_\kappa^m,H^1(A,\Q_l)) = Tr(F_\kappa^m,H^1(B,\Q_l))$ for all $m$.
Hence $A$ is isogenous to $B$  by the Tate conjectures for abelian varieties over finite fields. If $A$ and $B$ are isogenous, for the proof of corollary \ref{lc} over finite fields using our reduction to products of two simple abelian varieties, we can assume 
$A=B$ and the abstract subgroup $\bigcup_{n,m\in \Z} \{(na,ma) \ a\in A\}$ is Zariski dense in $X=A\times A$. So for every proper closed subset $D\subset X$ there exists $n,m$ with $(n,m)=1$ such that
$\varphi_M(A\times \{ 0\}) \cap D$ is a proper closed subset of $\varphi_M(A\times \{0\})$ for the automorphism $\varphi_M \in Aut(X)$ defined by a unimodular matrix $M\in Gl(2,\Z)$ with entries
$M_{11}=n$ and $M_{21}=m$. If the support of $K$ is $X$, then $supp({\cal H}^{-2\dim(A)}(K))=X$ and ${\cal H}^{-2\dim(A)}(K)$ is a smooth nonvanishing sheaf on $U=X\setminus Y$ for some
Zariski closed subset $Y\subset D$ with $D$ of dimension $d-1$. The intersection $U \cap \varphi_M(A\times\{0\})$
is Zariski dense in $\varphi_M(A\times\{0\})$ unless $\varphi_M(A\times\{0\}) \subset Y$. 
Hence without restriction of generality  $M^d= {}^p H^{-\dim(A)}(K\vert_{A\times \{0\}})$ can be assumed to have support  $A\times \{0\}$ by replacing $K$ with $\varphi_M^*(K)$ for some $M$. 

\medskip
This however contradicts step 15 of our last proof. In the characteristic zero case this proves the analog of corollary \ref{lc} below, in view of proposition \ref{Rec}. Another argument in characteristic zero, relying on proposition \ref{product}, was already given in proposition \ref{ToDo}.  Besides the assumptions (A1) and (A2), this is the only place where in the proof of the main result theorem \ref{MAINTHEOREM} of this paper the characteristic zero assumption enters, namely by using the support condition $supp(K)=X$ in the proof of proposition \ref{ToDo}. That this
support condition holds in characteristic zero followed from theorem \ref{Support} and its 
corollary \ref{Rec}. 

\medskip
Let us therefore give a third argument via Lang torsors, which completely avoids the use of the support condition $supp(K)=X$ made in proposition \ref{ToDo} and therefore, in view of the results of loc. cit.,  implies corollary \ref{lc} over finite fields:

\bigskip
Let $K_0$ respectively $K$ be perverse sheaves as in proposition \ref{product}.
We want to apply this proposition to the pullback $\wp^{(m)*}_0(K_0)$ of $K_0$ with respect to
some Lang morphism $\wp^{(m)}_0: X_0\times_{Spec(\kappa)}Spec(\kappa_m) \to X_0\times_{Spec(\kappa)} Spec(\kappa_m)$ respectively the corresponding isogeny over $k$
$$\xymatrix{ 0\ \ar[r] & X(\kappa_m) \ar[r] &   X   \ \ar[r]^{\wp^{(m)}} & \  X \ar[r] & \  0 } \ .$$ For this replace
$\kappa$ by $\kappa_m$ and $K$ by $\tilde K = \wp^{(m)}_*(K)$. Notice  
$\tilde K$ is defined over $\kappa_m$, and we write $\tilde K_0$ for the corresponding 
perverse sheaf on $X_0 \times_{\kappa} \kappa_m$. 
By [W2, cor. 7] the perverse sheaf $\tilde K$ is a simple monoidal perverse sheaf  with the same support as $K$ so that ${\cal H}^{-d}(\tilde K) =\delta_0$. So, once we replace $X_0,A_0,B_0$
by their scalar extension to $\kappa_m$ the perverse sheaf $\tilde K$ satisfies all the conditions
required for proposition \ref{product}. Since $\wp^{(m)}_*$ is a tensor functor, the properties $\tilde K^\vee \cong \tilde K$
and $\tilde K * \tilde K \cong H^\bullet \cdot \tilde K \oplus \tilde T$ hold for the translation invariant complex $\tilde T = \wp^{(m)}_*(T)$ on $B$, since the corresponding property holds for $K$. Finally $$ Rp_*(\tilde K) \ = \ \wp^{(m)}_* (Rp_*(K)) \ ,$$
because $p$ is defined over $\kappa$ and therefore commutes with $\wp^{(m)}$. 
Similarly, for characters $\varphi = \chi \circ \pi_1(\wp^{(m)})$, we get 
$ Rp_*(\tilde K_\chi)   =  \wp^{(m)}_* (Rp_*(K_{\varphi}))$.
Hence $Rp_*(\tilde K_{\chi_1}) $ is a translation invariant perverse sheaf on $B$ for almost all
characters $\chi_1$, since it is the direct image under $\wp^{(m)}_*$ of the corresponding translation invariant perverse sheaf $Rp_*(K_{\varphi})$ for $\varphi=(\varphi_1,\varphi_2)$ and $\varphi_1\notin \Sigma$. Therefore proposition \ref{product}
can be applied for $\tilde K$ instead of $K$ and again shows
$$  rank(Rp_*(\tilde K)) \ = \ 1\ .$$
Of course, also $rank(Rp_*(K))=1$ holds by proposition \ref{product}. 
Twisting $K$ by a suitable character, we may assume
that $Rp_*(K)$ is perverse translation-invariant of rank $r=1$ so that
$Rp_*(K) \cong \delta_B^\psi$. Hence we get 
$$ rank(Rp_*(\tilde K))  \ =\  rank(\wp^{(m)}_* (Rp_*(K))) \ =\  deg(\wp^{(m)}) \cdot rank(Rp_*(K)) \ =\  deg(\wp^{(m)}) \ $$
from $\wp^{(m)}_* (\delta_B^\psi)) = \oplus_{i=1}^{deg(\wp^{(m)})} (\delta_B^{\psi})_{\chi_i}$. 
This implies $ \# X(\kappa_m)= deg(\wp^{(m)})  = 1$, which is absurd for large $m$.
This contradiction implies that perverse sheaves $K$ with the properties required
as in proposition \ref{product} can not exist. Together with the other arguments used for the analogous proof in the characteristic zero case,  this implies

\begin{Corollary}\label{lc}
Let $k$ be the algebraic closure of a finite field $\kappa$. Suppose
any monoidal perverse Weil sheaf on a simple abelian variety over $k$ is translation invariant. Furthermore suppose that all simple constituents of an acyclic perverse sheaf on a simple abelian variety are acyclic.  Then 
any irreducible perverse Weil sheaf with Euler characteristic zero on an arbitrary abelian variety $X$ over $k$ is translation invariant under some abelian subvariety $Y \subseteq X$ of positive dimension.
\end{Corollary}

\bigskip\noindent

\medskip

\medskip

\bigskip\noindent
{\bf References}

\bigskip\noindent
[A] Abramovich D., {\it Subvarieties of semiabelian varieties}, Comp. Math. 90 (1994), 37 - 52

\medskip\noindent
[BBD] Beilinson A., Bernstein J., Deligne P., {\it Faiscaux pervers}, Asterisque 100 (1982) 

\medskip\noindent
[BK] B\"ockle G., Khare C., {\it Mod $l$ representations of arithmetic fundamental groups. II.
A conjecture of A.J de Jong}. Compos. Math. 142 (2006) 271 - 294

\medskip\noindent
[D] Drinfeld V., {\it On a conjecture of Kashiwara}, Mathematical Research Letters 8 (2001), 713 - 728

\medskip\noindent
[Dr2] Drinfeld V., {\it On a conjecture of Deligne}, Moscow Mathematical Journal 12 (2012)

\medskip\noindent
[Fin] Deligne P., {\it Theoremes de finitude en cohomologie $l$-adique}, in Seminaire de Geometrie Algebrique du Bois-Marie SGA 4$\frac{1}{2}$, {\it Cohomologie Etale}, Lecture Notes in Mathematics 569, Springer (1977)

\medskip\noindent
[FK] Franecki J., Kapranov M., {\it The Gauss map and a noncompact Riemann-Roch formula for constructible sheaves on semiabelian varieties}, arXiv:math/9909088

\medskip\noindent
[Gi] Ginsburg V., {\it Characteristic Varieties and vanishing cycles}, Invent. math. 84, 327 - 402 (1986) 

\medskip\noindent
[G] Gaitsgory D., {\it On de Jong's conjecture}, Israel J. Math. 157, no.1 (2007), 155 - 191

\medskip\noindent
[H] Hotta R., Takeuchi K., Tanisaki T., {\it D-Modules, Perverse sheaves, and representation theory}, Birkh\"auser Verlag (2008)

%
\bigskip\noindent
[KS] Kashiwara M., Schapira P., {\it Sheaves on manifolds}, Grundlehren der mathematischen Wissenschaften 292, Springer 2002

%

\bigskip\noindent
[KrW] Kr\"amer T., Weissauer R., {\it Vanishing theorems for constructible sheaves on abelian varieties}, arXiv:1111.4947v3 (2012)

\bigskip\noindent
[KW] Kiehl R., Weissauer R., {\it Weil conjectures, perverse sheaves and $l$-adic Fourier transform}, Erg. Math. Grenzg. (3.Folge) 42, Springer (2000)

\bigskip\noindent
[L] Lafforgue L., {\it Chtoucas de Drinfeld et correspondance de Langlands}, Invent. Math.
147 (2002), no. 1, 1 - 241

\bigskip\noindent
[M] Milne J.S. {\it Abelian Varieties}, in Arithmetic Geometry, edited by Cornell G., Silverman J.H., Springer (1986)

\bigskip\noindent
[BN] Weissauer R., {\it Brill-Noether sheaves}, arXiv:math/0610923v4 (2007)

\bigskip\noindent
[W] Weissauer R., {\it On Subvarieties of Abelian Varieties with degenerate Gau\ss\ mapping}, arXiv:1110.0095v2 (2011) 

\bigskip\noindent
[W1] Weissauer R., {\it A remark on the rigidity of BN-sheaves}, arXiv:1111.6095v1 (2011) 

\bigskip\noindent
[W2] Weissauer R., {\it On the rigidity of BN-sheaves}, arXiv:1204.1929 (2012) 

\bigskip\noindent
[W3] Weissauer R., {\it Why certain Tannaka groups attached to abelian varieties are almost connected}, preprint (2012)

\medskip\noindent

\medskip\noindent

\medskip\noindent

\end{document}